\documentstyle[a4, 12pt]{article}

\input amssym.def
\input amssym

\font\tenmsb=msbm10 scaled \magstep1
\font\sevenmsb=msbm7 scaled \magstep1
\font\fivemsb=msbm5 scaled \magstep1
\newfam\msbfam
\textfont\msbfam=\tenmsb
\scriptfont\msbfam=\sevenmsb
\scriptscriptfont\msbfam=\fivemsb
\def\Bbb#1{{\fam\msbfam\relax#1}}

\font\teneufm=eufm10 scaled \magstep1
\font\seveneufm=eufm7 scaled \magstep1
\font\fiveeufm=eufm5 scaled \magstep1
\newfam\eufmfam
\textfont\eufmfam=\teneufm
\scriptfont\eufmfam=\seveneufm
\scriptscriptfont\eufmfam=\fiveeufm
\def\frak#1{{\fam\eufmfam\relax#1}}

\newtheorem{proposition}{{\sc Proposition}}%[section]

\newtheorem{corollary}[proposition]{{\sc Corollary}}
\newtheorem{definition}[proposition]{{\sc Definition}}
\newtheorem{lemma}[proposition]{{\sc Lemma}}
\newtheorem{theorem}[proposition]{{\sc Theorem}}
\def\PROOF{{\noindent{\it Proof.\ }}}

\def\ident{\buildrel\sim\over\longrightarrow}

\def\build#1_#2^#3{{\mathrel{\mathop{\kern 0pt#1}\limits_{#2}^{#3}}}}
\def\posrel#1#2{\,{\build\hbox to 12mm{\rightarrowfill}_{#1}^{#2}}\,}
\def\hfld#1#2{\smash{\mathop{\hbox to 12mm{\rightarrowfill}}
     \limits^{\scriptstyle#1}_{\scriptstyle#2}}}
\def\hflg#1#2{\smash{\mathop{\hbox to 12mm{\leftarrowfill}}
     \limits^{\scriptstyle#1}_{\scriptstyle#2}}}

\def\diagramme#1{\def\normalbaselines{\baselineskip=0pt\lineskip=8pt
     \lineskiplimit=1pt}\matrix{#1}}

\title{Nearby cycles for local models of some Shimura varieties}
\author{T. Haines, B.C. Ng\^o}
\date{}
\begin{document}
\maketitle

\section{Introduction}

For certain classical groups $G$ and certain minuscule coweights $\mu$ of $G$, 
M. Rapoport and Th. Zink have constructed a projective scheme $M(G,\mu)$ over
${\Bbb Z}_p$ that is a local model for singularities at $p$ of some Shimura 
variety with level structure of Iwahori type at $p$. Locally for the \'{e}tale 
topology, $M(G,\mu)$ is isomorphic to a natural ${\Bbb Z}_p$-model 
${\mathcal M}(G,\mu)$ of the Shimura variety. 

The semi-simple trace of the Frobenius endomorphism on the nearby 
cycles of ${\mathcal M}(G,\mu)$ plays an important role in the computation 
of the local factor at $p$ of the semi-simple Hasse-Weil zeta function
of the Shimura variety, see \cite{Rapoport}. 
We can recover 
the semi-simple trace of Frobenius on the nearby cycles of 
${\mathcal M}(G,\mu)$ from that 
of the local model $M(G,\mu)$, see loc.cit. Thus the problem to calculate  
the function
$$x\in M(G,\mu)({\Bbb F}_q)\mapsto{\mathrm Tr}^{ss}({\mathrm Fr}_q,{\mathrm R}\Psi({{\bar{\Bbb Q}}_\ell})_x)$$
comes naturally.  R. Kottwitz has conjectured an explicit formula for this function.

To state this conjecture, we note that 
the set of ${\Bbb F}_q$-points of $M(G,\mu)$ can be naturally embedded as a finite set of Iwahori-orbits in the affine flag variety of $G({\Bbb F}_q(\!(t)\!))$ 
$$M(G,\mu)({\Bbb F}_q)\subset G\bigl({\Bbb F}_q(\!(t)\!)\bigr)/I$$
where $I$ is the standard Iwahori subgroup of 
$G\bigl({\Bbb F}_q(\!(t)\!)\bigr)$.

\vspace{5 pt}

\noindent {\bf Conjecture }({\bf Kottwitz}) {\em For all  $x\in M(G,\mu)({\Bbb F}_q)$},
$${\mathrm Tr}^{ss}({\mathrm Fr}_q,{\mathrm R}\Psi({{\bar{\Bbb Q}}_\ell})_x)=q^{\langle\rho,\mu\rangle }z_\mu(x).$$
\vspace{5pt}

\noindent Here $q^{\langle\rho,\mu\rangle }z_\mu(x)$ is the unique 
function in the center of the Iwahori-Hecke algebra of $I$-bi-invariant
functions with compact support in $G\bigl({\Bbb F}_p(\!(t)\!)\bigr)$, 
characterized by
$$q^{\langle\rho,\mu\rangle }z_\mu(x)*{\Bbb I}_K={\Bbb I}_{K\mu K}.$$
Here $K$ denotes the maximal compact subgroup $G({\Bbb F}_q[[t]])$ and 
${\Bbb I}_{K\mu K}$ denotes the characteristic function of the 
double-coset corresponding to a coweight $\mu$. 

Kottwitz' conjecture was first proved for the local model of a special type of Shimura variety with Iwahori type reduction at $p$ attached to the group ${\mathrm GL}(d)$ and minuscule coweight $(1,0^{d-1})$ (the ``Drinfeld case'') in \cite{Haines2}.  The method of that paper was one of direct computation: Rapoport had computed the function ${\mathrm Tr}^{ss}({\mathrm Fr}_q,{\mathrm R}\Psi({\bar {\Bbb Q}_\ell})_x)$ for the Drinfeld case (see \cite{Rapoport}), and so the result followed from a comparison with an explicit formula for the Bernstein function $z_{(1,0^{d-1})}$.  More generally, the explicit formula for $z_\mu$ in \cite{Haines2} is valid for any minuscule coweight $\mu$ of any quasisplit $p$-adic group.  
Making use of this formula, U. G\"{o}rtz verified Kottwitz' conjecture for a similar Iwahori-type Shimura variety attached to $G = {\mathrm GL}(4)$ and $\mu = (1,1,0,0)$, by computing the function ${\mathrm Tr}^{ss}({\mathrm Fr}_q,{\mathrm R}\Psi({\bar {\Bbb Q}_\ell})_x)$ 
for $x$ ranging over all 33 strata of the corresponding local model $M(G,\mu)$.

Shortly thereafter, A. Beilinson and D. Gaitsgory were motivated by Kottwitz' conjecture to attempt to produce all elements in the center of the Iwahori-Hecke algebra geometrically, via a nearby cycle construction.  For this they used Beilinson's deformation of the affine Grassmanian: a space over a curve $X$ whose fiber over a fixed point $x \in X$ is the affine flag variety of the group $G$, and whose fiber over every other point of $X$ is the affine Grassmanian of $G$.  In \cite{Gaitsgory} Gaitsgory proved a key commutativity result (similar to our Proposition 21) which is valid for any split group $G$ and any dominant coweight, in the function field setting.  His result also implies that the semi-simple trace of Frobenius on nearby cycles (of a $K$-equivariant perverse sheaf on the affine Grassmanian) corresponds to a function in the center of the Iwahori-Hecke algebra of $G$. 

The purpose of this article is to give a proof of Kottwitz' conjecture for the cases $G = {\mathrm GL}(d)$ and $G = {\mathrm GSp}(2d)$.  In fact we prove a stronger result (Theorem 11) which applies to arbitrary coweights, and which was also conjectured by Kottwitz (although only the case of minuscule coweights seems to be directly related to Shimura varieties).  

\vspace{5pt}
\noindent {\bf Main Theorem} {\em Let $G$ be either ${\mathrm GL}(d)$ or 
${\mathrm GSp}(2d)$.  Then for any dominant coweight $\mu$ of $G$, we have}
$$
{\mathrm Tr}^{ss}({\mathrm Fr}_q,{\mathrm R}\Psi^{M}({\mathcal A}_{\mu,\eta})) = (-1)^{2\langle\rho,\mu\rangle}\sum_{\lambda \leq \mu}m_\mu(\lambda)z_{\lambda}.
$$
\vspace{5pt}

Here $M$ is a member of an increasing family of schemes $M_{n_\pm}$ 
which contains the local models of Rapoport-Zink; the generic fiber of $M$ can be embedded in the affine Grassmanian of $G$, and ${\mathcal A}_{\mu,\eta}$ denotes the $K$-equivariant intersection complex corresponding to $\mu$.  The special fiber of $M$ embeds in the affine flag variety of 
$G({\bar{ {\Bbb F}}}_q(\!(t)\!))$ so we can think of the semi-simple trace of Frobenius on nearby cycles as a function in the Iwahori-Hecke algebra of $G$.

While the strategy of proof is similar to that of Beilinson and Gaitsgory, 
in order to get a statement which is valid over all local non-Archimedean fields we use a somewhat different model, based on spaces of lattices, in the construction of the schemes $M_{n_\pm}$ (we have not determined the precise relation between our model and that of Beilinson-Gaitsgory). This is necessary to compensate for the lack of an adequate notion of affine Grassmanian over $p$-adic fields. The union of the schemes $M_{n_\pm}$ can be thought of as a $p$-adic analogue of Beilinson's deformation of the affine Grassmanian.

We would like to thank G. Laumon and M. Rapoport for generous advice and 
encouragement. We would like to thank A. Genestier who has pointed 
out to us a mistake occurring in the first version of this paper.  We thank Robert Kottwitz for explaining the argument of Beilinson and Gaitsgory to us and for helpful conversations about this material. 

T. Haines acknowledges the hospitality and support of the Institut des Hautes \'{E}tudes Scientifiques in Bures-sur-Yvette in the spring of 1999, when this work was begun.  He is partially supported by an NSF-NATO Postdoctoral fellowship, and an NSERC Research grant.
 
B.C. Ng\^o was visiting the Max Planck Institut fuer Mathematik during the preparation of this article. 

\section{Rapoport-Zink local models}
\subsection{Some definitions in the linear case}
Let $F$ be a local non-Archimedean field. 
Let ${\mathcal O}$ denote the ring of integers 
of $F$ and let $k={\Bbb F}_q$ denote the residue field of ${\mathcal O}$. We choose a 
uniformizer $\varpi$ of ${\mathcal O}$. We denote by $\eta$ the generic point of 
$S={\mathrm Spec\,}({\mathcal O})$ and by $s$ its closed point.

For $G={\mathrm GL}(d)$ and for $\mu$ the minuscule coweight
$$(\underbrace{1,\ldots,1}_r,\underbrace{0,\ldots,0}_{d-r})$$
with $1\leq r\leq d-1$, the local model $M_{\mu}$ represents the functor
which associates to each ${\mathcal O}$-algebra $R$ the set of 
$L_\bullet=(L_0,\ldots,L_{d-1})$ where $L_0,\ldots,L_{d-1}$ are 
$R$-submodules of $R^d$ satisfying the following properties
\begin{itemize}
\item $L_0,\ldots,L_{d-1}$ are locally direct factors of corank $r$ in $R^d$, 
\item $\alpha'(L_0)\subset L_1,\,\alpha'(L_1)\subset L_2,\ldots,\,
\alpha'(L_{d-1})\subset L_0$ where $\alpha$ is the matrix
$$\alpha'=\pmatrix{0& 1&        &  \cr
               & \ddots& \ddots &  \cr
                    &  &   0    & 1\cr
                 \varpi&  &        & 0}$$
\end{itemize}
The projective $S$-scheme $M_{\mu}$ is a local model for singularities at $p$ 
of some Shimura variety for unitary group with level structure of 
Iwahori type at $p$ (see \cite{Rapoport},\cite{Rapoport-Zink}).    

Following a suggestion of G. Laumon, we introduce a new variable $t$ and 
rewrite the moduli problem of $M_{\mu}$ as follows.
Let  $M_{\mu}(R)$ be the set of $L_{\bullet}=(L_0,\ldots,L_{d-1})$
where $L_0,\ldots,L_{d-1}$ are 
$R[t]$-submodules of $R[t]^d/tR[t]^d$ satisfying the following properties
\begin{itemize}
\item as $R$-modules, $L_0,\ldots,L_{d-1}$ are locally direct factors of 
corank $r$ in $R[t]^d/tR[t]^d$, 
\item $\alpha(L_0)\subset L_1,\,\alpha(L_1)\subset L_2,\ldots,\,
\alpha(L_{d-1})\subset L_0$ where $\alpha$ is the matrix
$$\alpha=\pmatrix{0& 1&        &  \cr
              & \ddots& \ddots &  \cr
                   &  &   0    & 1\cr
              t+\varpi&  &        & 0}$$
\end{itemize}

Obviously, these two descriptions are equivalent because $t$ acts as $0$
on the quotient $R[t]^d/tR[t]^d$. 
Nonetheless, the latter description indicates how to construct larger $S$-schemes $M_\mu$, where $\mu$ runs over a certain cofinal family of dominant (nonminuscule) coweights. 

Let $n_-\leq 0 <n_+$ be two integers.

\begin{definition}
 Let $M_{r,n_{\pm}}$ be  the functor which
associates each ${\mathcal O}$-algebra $R$ the set of $L_{\bullet}=(L_0,\ldots,L_{d-1})$
where $L_0,\ldots,L_{d-1}$ are 
$R[t]$-submodules of $$t^{n_-}R[t]^d/t^{n_+}R[t]^d$$ 
satisfying the following properties
\begin{itemize}
\item as $R$-modules, $L_0,\ldots,L_{d-1}$ are locally direct factors with
rank $n_{+}d-r$ in $t^{n_-}R[t]^d/t^{n_+}R[t]^d$, 
\item $\alpha(L_0)\subset L_1,\,\alpha(L_1)\subset L_2,\ldots,\,
\alpha(L_{d-1})\subset L_0$.
\end{itemize}
\end{definition}

This functor is obviously represented by a closed sub-scheme in a product
of Grassmannians. In particular, $M_{r,n_{\pm}}$ is projective over $S$. 

In some cases, it is more convenient to adopt the following
equivalent description of the functor $M_{r,n_\pm}$. 
Let us consider $\alpha$ as an element of the group
$$\alpha\in{\mathrm GL}(d,{\mathcal O}[t,t^{-1},(t+\varpi)^{-1}]).$$
Let ${\mathcal V}_0,{\mathcal V}_1,\ldots,{\mathcal V}_d$ be the fixed ${\mathcal O}[t]$-submodules 
of ${\mathcal O}[t,t^{-1},(t+\varpi)^{-1}]^d$ defined by 
$${\mathcal V}_i=\alpha^{-i} {\mathcal O}[t]^d.$$ 
In particular, we have ${\mathcal V}_d=(t+\varpi)^{-1}{\mathcal V}_0$. Denote by ${\mathcal V}_{i,R}$ 
the tensor ${\mathcal V}_i\otimes_{\mathcal O} R$ for any ${\mathcal O}$-algebra $R$.

\begin{definition}
Let $M_{r,n_\pm}$ be the functor which associates to each ${\mathcal O}$-algebra
$R$ the set of 
$${\mathcal L}_\bullet=({\mathcal L}_0\subset{\mathcal L}_1\subset\cdots\subset
{\mathcal L}_d=(t+\varpi)^{-1}{\mathcal L}_0)$$
where ${\mathcal L}_0,{\mathcal L}_1,\ldots$ are $R[t]$-submodules of $R[t,t^{-1},(t+\varpi)^{-1}]^d$
satisfying the following conditions
\begin{itemize}
\item for all $i=0,\ldots,d-1$, we have 
$t^{n_+}{\mathcal V}_{i,R}\subset{\mathcal L}_i\subset t^{n_-}{\mathcal V}_{i,R}$, 
\item as $R$-modules, ${\mathcal L}_i/t^{n_+}{\mathcal V}_{i,R}$ is locally a direct factor of
$t^{n_-}{\mathcal V}_{i,R}/t^{n_+}{\mathcal V}_{i,R}$ with rank $n_+d-r$. 
\end{itemize}
\end{definition}

By using the isomorphism
$$\alpha^i:t^{n_-}{\mathcal V}_{i,R}/t^{n_+}{\mathcal V}_{i,R}
\ident t^{n_-}R[t]^d/t^{n_+}R[t]^d$$
we can associate to each sequence $L_\bullet=(L_i)$ as in Definition 1
 of $M_{r,n_\pm}$, the sequence ${\mathcal L}_\bullet=({\mathcal L}_i)$ 
as in Definition 2, in such a way that
$$\alpha^i({\mathcal L}_i/t^{n_+}{\mathcal V}_{i,R})=L_i.$$
This correspondence is clearly bijective. 
Therefore, the two definitions of the 
functor $M_{r,n_\pm}$ are equivalent.

It will be more convenient to consider the disjoint union  $M_{n_\pm}$
of projective schemes
$M_{r,n_\pm}$ for all $r$ for which $M_{r,n_\pm}$ makes sense, namely
$$M_{n_\pm}=\coprod_{dn_-\leq r \leq dn_+}M_{r,n_\pm},$$
instead of each connected component $M_{r,n_\pm}$ individually.

\subsection{Group action}
Definition 2 permits us to define a natural group action 
on $M_{n_\pm}$. Every $R[t]$-module ${\mathcal L}_i$ as above is included in
$$t^{n_+}R[t]^d\subset {\mathcal L}_i\subset t^{n_-}(t+\varpi)^{-1}R[t]^d.$$
Let ${\bar{\mathcal L}}_i$ denote its image in the quotient
$${\bar{\mathcal V}}_{n_\pm,R}=t^{n_-}(t+\varpi)^{-1}R[t]^d/t^{n_+}R[t]^d.$$
Obviously, ${\mathcal L}_i$ is completely determined by ${\bar{\mathcal L}}_i$. 

Let ${\bar{\mathcal V}}_i$ denote the image of 
${\mathcal V}_i$ in ${\bar{\mathcal V}}_{n_\pm}$.
We can view ${\bar{\mathcal V}}_{n_\pm}$
as the free $R$-module $R^{(n_+ -n_-+1)d}$ 
equipped with the endomorphism $t$ and with the filtration
$${\bar{\mathcal V}}_\bullet = ({\bar{\mathcal V}}_0\subset{\bar{\mathcal V}}_1
\cdots\subset{\bar{\mathcal V}}_d=(t+\varpi)^{-1}{\bar{\mathcal V}}_0)$$ 
which is stabilized by $t$. 

We now consider the functor ${J}_{n_\pm}$ which associates to each 
${\mathcal O}$-algebra $R$ the group ${J}_{n_\pm}(R)$ of all $R[t]$-automorphisms
of ${\bar{\mathcal V}}_{n_\pm}$ fixing the filtration 
${\bar{\mathcal V}}_\bullet$. This functor is represented by a closed subgroup
of ${\mathrm GL}((n_+-n_-+1)d)$ over $S$ that
acts in the obvious way on $M_{n_\pm}$.

\begin{lemma}
The group scheme $J_{n_\pm}$ is smooth over $S$.
\end{lemma}

\PROOF
Consider the functor ${\mathcal J}_{n_\pm}$ which associates to each 
${\mathcal O}$-algebra $R$ the ring ${\mathcal J}_{n_\pm}(R)$ of all $R[t]$-endomorphisms
of ${\bar{\mathcal V}}_{n_\pm}$ stabilizing the filtration 
${\bar{\mathcal V}}_\bullet$. This functor is obviously represented by a closed 
sub-scheme of the $S$-scheme 
\linebreak ${\frak{gl}}((n_+-n_-+1)d)$ of square matrices with
rank $(n_+ -n_-+1)d$.

The natural morphism of functors $J_{n_\pm}\rightarrow{\mathcal J}_{n_\pm}$ is an 
open immersion. Thus it suffices to prove that ${\mathcal J}_{n_\pm}$ is smooth
over $S$.

Giving an element of ${\mathcal J}_{n_\pm}$ is equivalent to giving $d$ 
vectors $v_1,\ldots,v_d$ such that 
$v_i\in t^{n_-}{\bar{\mathcal V}}_i$.
This implies that ${\mathcal J}_{n_\pm}$ is isomorphic to a trivial vector
bundle over $S$ of rank 
$$\sum_{i=1}^d{\mathrm rk}_{\mathcal O}(t^{n_-}{\mathcal V}_i/t^{n_+}{\mathcal O}[t]^d)=
   d^2(n_+ -n_- +1)-(d-1)d/2.$$ 
This finishes the proof of the lemma. $\square$

\subsection{Description of the generic fibre}

For this purpose, we use Definition 1 of $M_{n_\pm}$.
Let $R$ be an $F$-algebra. The matrix $\alpha$ then is invertible 
as an element
$$\alpha\in{\mathrm GL}(d,R[t]/t^{n_+-n_-}R[t]),$$
the group of automorphisms of $t^{n_-}R[t]^d/t^{n_+}R[t]^d$.

Let $(L_0,\ldots,L_{d-1})$ be an element of $M_{n_\pm}(R)$.
As $R$-modules, the $L_i$ are locally direct factors of the same rank.
For $i=1,\ldots,d-1$, the inclusion
$\alpha(L_{i-1})\subset L_i$
implies the equality
$\alpha(L_{i-1})=L_i$.
In this case, the last inclusion
$\alpha(L_{d-1})\subset L_0$
is automatically an equality, because the matrix
$$\alpha^d={\mathrm diag\,}(t+\varpi,\ldots,t+\varpi)$$
satisfies the property: $\alpha^d(L_0)=L_0$.
In others words, the whole sequence $(L_0,\ldots,L_{d-1})$ is completely
determined by $L_0$. 

Let us reformulate the above statement in a more precise way.
Let ${\mathrm Grass}_{n_\pm}$ be the functor which associates to each 
${\mathcal O}$-algebra $R$ the set of $R[t]$-submodules $L$ of  
$t^{n_-}R[t]^d/ t^{n_+}R[t]^d$ which, as $R$-modules, are locally
direct factors of $t^{n_-}R[t]^d/ t^{n_+}R[t]^d$.
Obviously, this functor is represented by a closed sub-scheme of a 
Grassmannian. In particular, it is proper over $S$.

Let $\pi: M_{n_\pm}\rightarrow {\mathrm Grass}_{n_\pm}$
be the morphism defined by
$$\pi(L_0,\ldots,L_{d-1})=L_0.$$

The above discussion can be reformulated as follows.

\begin{lemma}
The morphism $\pi:M_{n_\pm}\rightarrow{\mathrm Grass}_{n_\pm}$ 
is an isomorphism over the generic point $\eta$ of $S$. $\square$
\end{lemma}

Let $K_{n_\pm}$ the functor which associates to each ${\mathcal O}$-algebra $R$ 
the group 
$$K_{n_\pm}={\mathrm GL}(d,R[t]/t^{n_+-n_-}R[t]).$$
Obviously, it is represented by a smooth group scheme over $S$
and acts naturally on ${\mathrm Grass}_{n_\pm}$.
This action yields a decomposition into orbits that are smooth over $S$
$${\mathrm Grass}_{n_\pm}=\coprod_{\lambda\in\Lambda(n_\pm)} O_\lambda$$
where $\Lambda(n_\pm)$ is the finite set
of sequences of integers $\lambda=(\lambda_1,\ldots,\lambda_d)$
satisfying the following condition
$$ n_+\geq\lambda_1\geq\cdots\geq\lambda_d\geq n_-.$$
This set $\Lambda(r,n_\pm)$ can be viewed as a finite subset of 
the cone of dominant coweights of $G={\mathrm GL}(d)$ and conversely, every 
dominant coweight of $G$ occurs in some $\Lambda(n_\pm)$.
For all $\lambda\in \Lambda(n_\pm)$, we have
$$O_\lambda(F)=K_F\,t^{\lambda} K_F/K_F.$$
Here $K_F={\mathrm GL}(d,F[[t]])$ is the standard maximal ``compact'' 
subgroup of $G_F={\mathrm GL}(d,F(\!(t)\!))$
and acts on ${\mathrm Grass}_{n_\pm}(F)$ through the quotient
$K_{n_\pm}(F)$. The above equality holds if one replaces $F$ by any field which is 
also an ${\mathcal O}$-algebra, since $K_{n_\pm}$ is smooth;
in particular it holds for the residue field $k$.  

We derive from  the above lemma the description
$$M_{n_\pm}(F)=\coprod_{\lambda\in\Lambda(n_\pm)}K_F\, t^{\lambda} K_F/K_F.$$ 

We will need to compare the action of $J_{n_\pm}$ on $M_{n_\pm}$
and the action of $K_{n_\pm}$ on ${\mathrm Grass}_{n_\pm}$. 
By definition, $J_{n_\pm}(R)$ is a subgroup of 
$$J_{n_\pm}(R)\subset{\mathrm GL}(d,R[t]/t^{n_+-n_-}(t+\varpi)R[t])$$
for any ${\mathcal O}$-algebra $R$.
By using the natural homomorphism
$${\mathrm GL}(d,R[t]/t^{n_+-n_-}(t+\varpi)R[t])\rightarrow
{\mathrm GL}(d,R[t]/t^{n_+-n_-}R[t])$$
we get a homomorphism $J_{n_\pm}(R)\rightarrow K_{n_\pm}(R)$.
This gives rises to a homomorphism of group schemes 
$\rho:J_{n_\pm}\rightarrow K_{n_\pm}$,
which is surjective over the generic point $\eta$
of $S$.

The proof of the following lemma is straightforward.

\begin{lemma}
With respect to the homomorphism $\rho:J_{n_\pm}\rightarrow K_{n_\pm}$,
and to the morphism $\pi:M_{n_\pm}\rightarrow {\mathrm Grass}_{n_\pm}$,
the action  of $J_{n_\pm}$ on $M_{n_\pm}$
and the action of $K_{n_\pm}$ on ${\mathrm Grass}_{n_\pm}$ are compatible. $\square$
\end{lemma}

\subsection{Description of the special fibre}

For this purpose, we will use Definition 2 of $M_{n_\pm}$.
The functor $M_{r, n_\pm}$ associates to each $k$-algebra $R$ the set of
$${\mathcal L}_\bullet=({\mathcal L}_0\subset{\mathcal L}_1\subset\cdots\subset
{\mathcal L}_d=t^{-1}{\mathcal L}_0)$$
where ${\mathcal L}_0,{\mathcal L}_1,\ldots$ are $R[t]$-submodules of $R[t,t^{-1}]^d$
satisfying the following conditions
\begin{itemize}
\item for all $i=0,\ldots,d-1$, we have 
$t^{n_+}{\mathcal V}_{i,R}\subset{\mathcal L}_i\subset t^{n_-}{\mathcal V}_{i,R}$, 
\item as an $R$-module, each ${\mathcal L}_i/t^{n_+}{\mathcal V}_{i,R}$ is locally a direct factor of
$t^{n_-}{\mathcal V}_{i,R}/t^{n_+}{\mathcal V}_{i,R}$ with rank $n_+d-r$. 
\end{itemize}

Let $I_{k}$ denote the standard Iwahori subgroup of 
$G_k={\mathrm GL}\bigl(d,k(\!(t)\!)\bigr)$, that is, 
the subgroup of integer matrices ${\mathrm GL}\bigl(d,k[[t]]\bigr)$ whose
reduction mod $t$ lies in the subgroup of upper triangular matrices in 
${\mathrm GL}(d,k)$. The set of $k$-points of $M_{n_\pm}$ can be realized as a finite 
subset in the set of affine flags of ${\mathrm GL}(d)$
$$M_{n_\pm}(k)\subset G_k/I_k.$$ 
By definition, the $k$-points of $J_{n_\pm}$ are the matrices in
${\mathrm GL}(d,k[t]/t^{n_+-n_-+1}k[t])$ whose reduction mod $t$ is upper
triangular. Thus, $J_{n_\pm}(k)$ is a quotient of $I_k$.
Obviously, the action of $J_{n_\pm}(k)$ on $M_{n_\pm}(k)$
and the action of $I_k$ on $G_k/I_k$ are compatible. 
Therefore, for each $r$ such that $dn_- \leq r \leq dn_+$ there exists a finite subset 
${\tilde W}(r, n_\pm)\subset {\tilde W}$ of the affine Weyl 
group ${\tilde W}$ such that
$$M_{n_\pm}(k)=\coprod_{w\in {\tilde W}(n_\pm)} I_k wI_k/I_k,$$
where ${\tilde W}(n_\pm) = \coprod_{r} {\tilde W}(r, n_\pm)$.
One can see easily that any element $w\in{\tilde W}$ occurs in the
finite subset ${\tilde W}(n_\pm)$ for some $n_\pm$.  But 
the exact determination of the finite sets ${\tilde W}(r,n_\pm)$ is 
a difficult combinatorial problem; 
for the case of minuscule coweights of ${\mathrm GL}(d)$ (i.e., $n_+ = 1$ and $n_- = 0$) 
these sets have been described by Kottwitz and Rapoport (\cite{Kora}).

Let us recall that
$${\mathrm Grass}_{n_\pm}(k)=\coprod_{\lambda\in\Lambda(n_\pm)} 
K_k\,t^{\lambda} K_k/K_k.$$
The proof of the next lemma is straightforward.

\begin{lemma}
The map $\pi(k):M_{n_\pm}(k)\rightarrow{\mathrm Grass}_{n_\pm}(k)$ is the restriction
of the natural map $G_k/I_k\rightarrow G_k/K_k$.
\end{lemma}

\subsection{Symplectic case}

For the symplectic case, we will give only the definitions of 
the symplectic analogues of the objects which were considered 
in the linear case. The statements of Lemmas 3,4,5 and 
6 remain unchanged.

In this section, the group $G$ stands for ${\mathrm GSp}(2d)$ 
associated to the symplectic form $\langle\, ,\, \rangle $ represented by the matrix
$$\pmatrix{0 & J\cr -J & 0}$$
where $J$ is the anti-diagonal matrix with entries equal to $1$. 
Let $\mu$ denote  the minuscule coweight 
$$\mu=(\underbrace{1,\ldots,1}_{d}, \underbrace{0,\ldots,0}_{d}).$$

Following Rapoport and Zink (\cite{Rapoport-Zink})
the local model $M_\mu$ represents the functor which associates
to each ${\mathcal O}$-algebra $R$ the set of sequences $L_\bullet=(L_0,\ldots,L_d)$
where $L_0,\ldots,L_d$ are $R$-submodules of $R^{2d}$ 
satisfying the following properties

\begin{itemize}
\item $L_0,\ldots,L_d$ are locally direct factors of $R^{2d}$ of rank $d$,
\item $\alpha'(L_0)\subset L_1,\ldots,\alpha'(L_{d-1})\subset L_d$ where 
$\alpha'$ is the matrix of size $2d\times 2d$
$$\alpha'=\pmatrix{0& 1&        &  \cr
               & \ddots& \ddots &  \cr
                    &  &   0    & 1\cr
                 \varpi&  &        & 0}$$
\item $L_0$ and $L_d$ are isotropic with respect to $\langle\, ,\, \rangle $.
\end{itemize}

Just as in the linear case, let us introduce a new variable $t$
and give the symplectic analog of Definition 2.
We consider the matrix of size $2d\times 2d$ 
$$\alpha'=\pmatrix{0& 1&        &  \cr
               & \ddots& \ddots &  \cr
                    &  &   0    & 1\cr
                t+\varpi&  &        & 0}$$
viewed as an element of
$$\alpha\in{\mathrm GL}(2d,{\mathcal O}[t,t^{-1},(t+\varpi)^{-1}]).$$
Denote by ${\mathcal V}_0,\ldots,{\mathcal V}_{2d-1}$ the fixed ${\mathcal O}[t]$-submodules 
of ${\mathcal O}[t,t^{-1},(t+\varpi)^{-1}]^{2d}$ defined by
${\mathcal V}_i=\alpha^{-i}{\mathcal O}[t]^{2d}$. For an ${\mathcal O}$-algebra $R$, 
let ${\mathcal V}_{i,R}$ denote ${\mathcal V}_i\otimes_{\mathcal O} R$.

For any $R[t]$-submodule ${\mathcal L}$ of $R[t,t^{-1},(t+\varpi)^{-1}]^{2d}$,
the $R[t]$-module
$${\mathcal L}^{\perp'} = \{x\in R[t,t^{-1},(t+\varpi)^{-1}]^{2d}\mid
\forall y\in {\mathcal L}, t^n(t+\varpi)^{n'}\langle x,y\rangle \in R[t]\}$$
is called the {\em dual} of ${\mathcal L}$ with respect to the form
$\langle\,,\,\rangle ' = t^n(t+\varpi)^{n'}\langle\,,\,\rangle $.
Thus ${\mathcal V}_0$ is autodual with respect to the form 
$\langle\,,\,\rangle $ and ${\mathcal V}_d$ is  autodual with respect to the form 
$(t+\varpi)\langle\,,\,\rangle $.

Here is the symplectic analog of Definition 2 of the model $M_{n_\pm}$.  
For $n_-=0$ and $n_+=1$, $M_{n_\pm}$ will coincide with $M_\mu$, for $\mu = (1^d, 0^d)$:

\begin{definition} For any $n_-\leq 0<n_+$,
let $M_{n_\pm}$ be the functor which associates to each ${\mathcal O}$-algebra
$R$ the set of sequences
$${\mathcal L}_\bullet=({\mathcal L}_0\subset{\mathcal L}_1\subset\cdots\subset{\mathcal L}_d)$$
where ${\mathcal L}_0,\ldots,{\mathcal L}_d$ are $R[t]$-submodules of 
$R[t,t^{-1},(t+\varpi)^{-1}]^{2d}$ satisfying the following properties
\begin{itemize}
\item for all $i=0,\ldots,d$, we have 
$t^{n_+}{\mathcal V}_{i,R}\subset{\mathcal L}_i\subset t^{n_-}{\mathcal V}_{i,R}$,
\item as $R$-modules, ${\mathcal L}_i/t^{n_+}{\mathcal V}_{i,R}$ is locally a direct factor of
$t^{n_-}{\mathcal V}_{i,R}/t^{n_+}{\mathcal V}_{i,R}$ of rank $(n_+-n_-)d$,
\item ${\mathcal L}_0$  is autodual with respect to the form
$t^{-n_--n_+}\langle\, ,\, \rangle $, and 
${\mathcal L}_d$ is autodual with respect to the form 
$t^{-n_--n_+}(t+\varpi)\langle\, ,\, \rangle $.
\end{itemize}
\end{definition}

Let us now define the natural group action on $M_{n_\pm}$.
The functor $J_{n_\pm}$ associates to each ${\mathcal O}$-algebra
$R$ the group $J_{n_\pm}(R)$ of $R[t]$-linear automorphisms of  
$${\bar{\mathcal V}}_{n_\pm,R}=t^{n_-}(t+\varpi)^{-1}R[t]^{2d}/t^{n_+}R[t]^{2d}$$
which fix the filtration
$${\bar{\mathcal V}}_{\bullet,R}=({\bar{\mathcal V}}_{0,R}\subset\cdots\subset{\bar{\mathcal V}}_{d,R})$$
(the image of ${\mathcal V}_{\bullet,R}$ in ${\bar {\mathcal V}}_{n_\pm,R}$) and which fix the symplectic form
$t^{-n_--n_+}(t+\varpi)\langle\,,\,\rangle $, up to a unit in $R$.
This functor is represented by an $S$-group scheme $J_{n_\pm}$
which acts on $M_{n_\pm}$.  Lemma 3 remains true in the symplectic case :
$J_{n_\pm}$ is a {\em smooth} group scheme over $S$. The proof
is completely similar to the linear case.

\bigskip
Let us now describe the generic fibre of $M_{n_\pm}$.
Let ${\mathrm Grass}_{n_\pm}$ be the functor which associates to 
each ${\mathcal O}$-algebra $R$
the set of $R[t]$-submodules $L$ of $t^{n_-}R[t]^{2d}/t^{n_+}R[t]^{2d}$
which, as $R$-modules, are locally direct factors of rank $(n_+-n_-)d$
and which are isotropic with respect to $t^{-n_--n_+}\langle\,,\,\rangle $.
Then the morphism $\pi:M_{n_\pm}\rightarrow{\mathrm Grass}_{n_\pm}$ 
defined by $\pi(L_\bullet)=L_{0}$ is an isomorphism over the generic
point $\eta$ of $S$.  Let $K_{n_\pm}$ denote the functor which associates to each 
${\mathcal O}$-algebra $R$ the group of 
$R[t]$-automorphisms of $t^{n_-}R[t]^{2d} / t^{n_+}R[t]^{2d}$ which fix the symplectic 
form $t^{-n_- -n_+}\langle\,,\,\rangle$ up to a unit in $R$.  Then $K_{n_\pm}$ is represented by a smooth group scheme over $S$, and it acts in the obvious way on 
${\mathrm Grass}_{n_\pm}$.  Consequently, we have a stratification in orbits
of the generic fibre $M_{n_\pm,\eta}$
$$M_{n_\pm,\eta}=\coprod_{\lambda\in \Lambda(n_\pm)}O_{\lambda,\eta}.$$
Here $\Lambda(n_\pm)$ is the set of sequences 
$\lambda=(\lambda_1,\ldots,\lambda_d)$ satisfying
$$n_+\geq\lambda_1\geq\cdots\geq\lambda_d\geq {n_++n_-\over 2},$$
and can be viewed as finite subset of the cone of dominant coweights
of $G={\mathrm GSp}(2d)$. One can easily check that each dominant coweight 
of ${\mathrm GSp}(2d)$ occurs in some $\Lambda(n_\pm)$.
For any $\lambda\in \Lambda(n_\pm)$, we have also
$$O_{\lambda,\eta}(F)=K_F t^\lambda K_F/K_F$$
where $K_F=G(F[[t]])$ is the ''maximal compact'' subgroup 
of $G_F=G(F(\!(t)\!))$. 

\bigskip
Next we turn to the special fiber of $M_{n_\pm}$.  For this it is most convenient to give a slight reformulation of Definition 7 above.  Let $R$ be any ${\mathcal O}$-algebra.  It is easy to see that specifying a sequence ${\mathcal L}_\bullet = ({\mathcal L}_0 \subset \dots {\mathcal L}_d)$ as in Definition 7 is the same as specifying a periodic ``lattice chain'' 
$$
\dots \subset {\mathcal L}_{-1} \subset {\mathcal L}_0 \subset \dots \subset {\mathcal L}_{2d} = (t+\varpi)^{-1}{\mathcal L}_0 \subset \dots
$$
consisting of $R[t]$-submodules of $R[t,t^{-1},(t+\varpi)^{-1}]^{2d}$ with the following properties:

\begin{itemize}
\item $t^{n_+}{\mathcal V}_{i, R} \subset {\mathcal L}_i \subset t^{n_-}{\mathcal V}_{i, R}$, where ${\mathcal V}_{i, R} = \alpha^{-i}{\mathcal V}_{0, R}$, for every $i \in {\Bbb Z}$,
\item ${\mathcal L}_i / t^{n_+}{\mathcal V}_{i, R}$ is locally a direct factor of rank $(n_+ - n_-)d$, for every $i \in {\Bbb Z}$,
\item ${\mathcal L}^{\perp}_i = t^{-n_- - n_+}{\mathcal L}_{-i}$, for every $i \in {\Bbb Z}$,
\end{itemize}
where $\perp$ is defined using the original symplectic form $\langle\,,\,\rangle$ on $R[t,t^{-1},(t+\varpi)^{-1}]^{2d}$.  
We denote by $I_k$ the standard Iwahori subgroup of ${\mathrm GSp}(2d,k[[t]])$, 
namely, the stabilizer in this group of the periodic lattice chain 
${\mathcal V}_{\bullet, k[[t]]}$.
There is a canonical surjection $I_k \rightarrow J_{n_\pm}(k)$ and so the Iwahori subgroup $I_k$ acts via its quotient $J_{n_\pm}(k)$ on the set $M_{n_\pm}(k)$.  Moreover, the $I_k$-orbits in $M_{n_\pm}(k)$ are parametrized by a certain finite set ${\tilde W}(n_\pm)$ of the affine Weyl group ${\tilde W}({\mathrm GSp}(2d))$
$$
M_{n_\pm}(k) = \coprod_{w \in {\tilde W}(n_\pm)} I_k \, w \, I_k / I_k.
$$
The precise description of the sets ${\tilde W}(n_\pm)$ is a difficult combinatorial problem (see \cite{Kora} for the case $n_+=1$, $n_-=0$), but one can easily see that 
any $w \in {\tilde W}({\mathrm GSp}(2d))$ is contained in some ${\tilde W}(n_\pm)$.

The definitions of the group scheme action of $K_{n_\pm}$ on ${\mathrm Grass}_{n_\pm}$, 
of the homomorphism $\rho:J_{n_\pm}\rightarrow K_{n_\pm}$
and the compatibility properties (Lemmas 5,6) are obvious and will be
left to the reader.

\section{Semi-simple trace on nearby cycles}

\subsection{Semi-simple trace}
The notion of semi-simple trace was introduced by Rapoport in \cite{Rapoport} and its good
properties were mentioned there. The purpose of this section is only to 
give a more systematic presentation in insisting on the important fact that
the semi-simple trace furnish a kind of sheaf-function dictionary \`a la
Grothendieck. In writing this section, we have benefited from very helpful
explanations of Laumon. 

Let ${\bar F}$ be a separable closure of the local field $F$.  
Let $\Gamma$ be the Galois group ${\mathrm Gal}({\bar F}/F)$ 
of $F$ and let $\Gamma_0$ be the inertia 
subgroup of $\Gamma$ defined by the exact sequence
$$1\rightarrow \Gamma_0\rightarrow\Gamma\rightarrow{\mathrm Gal}({\bar k}/k)\rightarrow 1.$$
For any prime $\ell\not= p$, there exists a canonical surjective homomorphism
$$t_\ell:\Gamma_0\rightarrow {\Bbb Z}_\ell(1).$$

Let ${\mathcal R}$ denote the abelian category of continuous, 
finite dimensional $\ell$-adic
representations of $\Gamma$. Let $(\rho,V)$ be an object of ${\mathcal R}$
$$\rho:\Gamma\rightarrow{\mathrm GL}(V).$$
According to Grothendieck, the restricted representation $\rho(\Gamma_0)$ is 
{\em quasi-unipotent} i.e. there exists a finite-index subgroup $\Gamma_1$ of $\Gamma_0$
which acts unipotently on $V$ (the residue field $k$ is supposed finite).
There exists then an unique nilpotent morphism, the {\em logarithm} of $\rho$
$$N:V\rightarrow V(-1)$$
characterized by the following property : for all $g\in \Gamma_1$, we have
$$\rho(g)=\exp(N t_\ell(g)).$$

Following Rapoport, an increasing filtration ${\mathcal F}$ of $V$ will be
called {\em admissible} if it is stable under the action of $\Gamma$ and
such that $\Gamma_0$ operates on the associated graded ${\mathrm gr}^{\mathcal F}_\bullet(V)$
through a finite quotient. Admissible filtrations always exist :
we can take for instance the filtration defined by the kernels of the powers 
of $N$.

We define the semi-simple trace of Frobenius on 
$V$ as
$${\mathrm Tr}^{ss}({\mathrm Fr}_q,V)=\sum_k {\mathrm Tr}({\mathrm Fr}_q,{\mathrm gr}^{\mathcal F}_k(V)^{\Gamma_0}) .$$

\begin{lemma}
The semi-simple trace ${\mathrm Tr}^{ss}({\mathrm Fr}_q,V)$ does not depend on the choice 
of the admissible filtration ${\mathcal F}$.
\end{lemma}

\PROOF Let us first consider the case where $\Gamma_0$ acts on $V$ through a 
finite quotient. Since the functor taking invariant of a finite group
acting on ${{\bar{\Bbb Q}}_\ell}$-vector space is exact, the graded associated to the filtration
${\mathcal F}'$ of $V^{\Gamma_0}$ induced by ${\mathcal F}$ is equal to 
${\mathrm gr}^{\mathcal F}_\bullet(V)^{\Gamma_0}$
$${\mathrm gr}_k^{\mathcal F'}(V^{\Gamma_0})={\mathrm gr}^{\mathcal F}_k(V)^{\Gamma_0}.$$
Consequently
$${\mathrm Tr}({\mathrm Fr}_q,V^{\Gamma_0})=\sum_k {\mathrm Tr}
({\mathrm Fr}_q,{\mathrm gr}^{\mathcal F}_k(V)^{\Gamma_0}) .$$

In the general case, any two admissible filtrations admit a third 
finer admissible filtration. By using the above case, one sees the semi-simple
trace associated to each of the two first admissible filtrations is 
equal to the semi-simple trace associated to the third one and
the lemma follows. $\square$

\begin{corollary} The function defined by
$$V\mapsto{\mathrm Tr}^{ss}({\mathrm Fr}_q,V)$$
on the set of isomorphism classes $V$ of ${\mathcal R}$,
factors through the Grothendieck group of ${\mathcal R}$.
\end{corollary}

For any object $C$ of the derived category associated to ${\mathcal R}$,
we put
$${\mathrm Tr}^{ss}({\mathrm Fr}_q,C)=\sum_i (-1)^i {\mathrm Tr}^{ss}({\mathrm Fr}_q,{\mathrm H}^i(C)).$$
By the above corollary, for any distinguished triangle 
$$C\rightarrow C'\rightarrow C''\rightarrow C[1]$$
the equality
$${\mathrm Tr}^{ss}({\mathrm Fr}_q,C)+{\mathrm Tr}^{ss}({\mathrm Fr}_q,C'')={\mathrm Tr}^{ss}({\mathrm Fr}_q,C')$$
holds.

Let $X$ be a $k$-scheme of finite type, $X_{\bar s}=X\otimes_k{\bar k}$.
Let $D^b_c(X\times_k\eta)$ the derived category associated to the abelian
category of constructible $\ell$-adic sheaves
on $X_{\bar s}$ equipped with an action of $\Gamma$ compatible with the
action of $\Gamma$ on $X_{\bar s}$ through ${\mathrm Gal}({\bar k}/k)$, 
see \cite{Deligne}. Let ${\mathcal C}$ be an object of 
$D^b_c(X\times_k\eta)$. For 
any $x\in X(k)$, the fibre ${\mathcal C}_x$ is an object of the derived
category of ${\mathcal R}$.
Thus we can define the function semi-simple trace
$$\tau^{ss}_{\mathcal C}:X(k)\rightarrow{{\bar{\Bbb Q}}_\ell}$$
by
$$\tau^{ss}_{\mathcal C}(x)={\mathrm Tr}^{ss}({\mathrm Fr}_q,{\mathcal C}_x).$$

This association ${\mathcal C}\mapsto \tau^{ss}_{\mathcal C}$ furnishes an analog 
of the usual sheaf-function dictionary of Grothendieck (see 
\cite{Grothendieck}):

\begin{proposition}
Let $f:X\rightarrow Y$ be a morphism between $k$-schemes of finite type.
\begin{enumerate}
\item Let ${\mathcal C}$ be an object of $D^b_c(Y\times_k\eta)$. For all 
$x\in X(k)$, we have
$$\tau^{ss}_{f^*{\mathcal C}}(x)=\tau^{ss}_{\mathcal C}(f(x))$$
\item Let ${\mathcal C}$ be an object of $D^b_c(X\times_k\eta)$. For all
$y\in Y(k)$, we have
$$\tau^{ss}_{{\mathrm R}f_!{\mathcal C}}(y)=\sum_{\scriptstyle x\in X(k)\atop
\scriptstyle f(x)=y} \tau^{ss}_{\mathcal C}(x).$$
\end{enumerate}
\end{proposition}

\PROOF
The first statement is obvious because $f^*{\mathcal C}_x$ and 
${\mathcal C}_{f(x)}$ are canonically isomorphic as objects of the derived 
category of ${\mathcal R}$.

It suffices to prove the second statement in the case $Y=s$. 
By Corollary 9 and ``d\'{e}vissage'', it 
suffices to consider the case where ${\mathcal C}$
is concentrated in only one degree, say in the degree zero. 
Denote $C={\mathcal H}^0({\mathcal C})$ and choose an admissible filtration of $C$
$$0=C_0\subset C_1\subset C_2\subset\cdots\subset C_n=C.$$

The associated spectral sequence
$${\mathrm E}^{i,j-i}_1={\mathrm H}^j_c(X_{\bar s}, C_i/C_{i-1})
\Longrightarrow {\mathrm H}^j_c(X_{\bar s},C)$$
yields an abutment filtration on ${\mathrm H}^j_c(X_{\bar s},C)$
with associated graded ${\mathrm E}_\infty^{i,j-i}$.
Since the inertia group acts on ${\mathrm E}^{i,j-i}_1$ through a finite quotient,
the same property holds for ${\mathrm E}_\infty^{i,j-i}$ because 
${\mathrm E}_\infty^{i,j-i}$ is a subquotient of ${\mathrm E}^{i,j-i}_1$.
Consequently, the abutment filtration on ${\mathrm H}^j_c(X_{\bar s},C)$
is an admissible filtration and by definition, we have
$${\mathrm Tr}^{ss}({\mathrm Fr}_q,Rf_!C)=\sum_{i,j}(-1)^j
{\mathrm Tr}({\mathrm Fr}_q,({\mathrm E}_\infty^{i,j-i})^{\Gamma_0}).$$

Now, the identity in the Grothendieck group
$$\sum_{i,j}(-1)^j{\mathrm E}_1^{i,j-i}=
\sum_{i,j}(-1)^j{\mathrm E}_\infty^{i,j-i}$$
implies
$$\sum_{i,j}(-1)^j({\mathrm E}_1^{i,j-i})^{\Gamma_0}
=\sum_{i,j}(-1)^j{(\mathrm E}_\infty^{i,j-i})^{\Gamma_0}$$
because taking the invariants by finite group
is an exact functor.

The same exactness implies
$${(\mathrm E}_1^{i,j-i})^{\Gamma_0}={\mathrm H}^j_c(X_{\bar s}, C_{i}/C_{i-1})^{\Gamma_0}
={\mathrm H}^j_c(X_{\bar s}, (C_{i}/C_{i-1})^{\Gamma_0}).$$
By putting the above equalities together, we obtain
$${\mathrm Tr}^{ss}({\mathrm Fr}_q,Rf_!C)=\sum_{i,j}(-1)^j
{\mathrm Tr}({\mathrm Fr}_q,{\mathrm H}^j_c(X_{\bar s}, (C_{i}/C_{i-1})^{\Gamma_0})).$$

By using now the Grothendieck-Lefschetz formula, we have
$$\sum_{x\in X(k)}{\mathrm Tr}({\mathrm Fr}_q,(C_{i}/C_{i-1})^{\Gamma_0}_x)
=\sum_{j}(-1)^j
{\mathrm Tr}({\mathrm Fr}_q,{\mathrm H}^j_c(X_{\bar s}, (C_{i}/C_{i-1})^{\Gamma_0})).$$
Consequently,
$${\mathrm Tr}^{ss}({\mathrm Fr}_q,Rf_!C)=\sum_{x\in X(k)}{\mathrm Tr}^{ss}({\mathrm Fr}_q,C_x).\ \ \square$$

\subsection{Nearby cycles}

Let ${\bar\eta}={\mathrm Spec\,}({\bar F})$ 
denote the geometric generic point of $S$,
${\bar S}$ be the normalization of $S$ in ${\bar\eta}$ and ${\bar s}$ be
the closed point of ${\bar S}$. For an $S$-scheme $X$ of finite type, 
let us denote by ${\bar\jmath}^X:X_{\bar\eta}\rightarrow X_{\bar S}$ 
the morphism deduced by base change from
${\bar\jmath}:{\bar\eta}\rightarrow {\bar S}$ and denote by 
${\bar\imath}^X:X_{\bar s}\rightarrow X_{\bar S}$ 
the morphism deduced from ${\bar\imath}:{\bar s}\rightarrow{\bar S}$.

The nearby cycles of an $\ell$-adic complex $C_\eta$ on $X_\eta$, 
is the complex of $\ell$-adic sheaves defined by
$${\mathrm R}\Psi^X(C_\eta)=
i^{X,*}{\mathrm R}{\bar\jmath}^X_{*} {\bar\jmath}^{X,*} C_\eta.$$
The complex ${\mathrm R}\Psi^X(C_\eta)$ is equipped with an action of 
$\Gamma$ compatible with the action of $\Gamma$ on $X_{\bar s}$
through the quotient ${\mathrm Gal}({\bar k}/k)$.

For $X$ a proper $S$-scheme, we have a canonical isomorphism
$${\mathrm R}\Gamma(X_{\bar s},{\mathrm R}\Psi(C_\eta))
={\mathrm R}\Gamma(X_{\bar\eta},C_\eta)$$
compatible with the natural actions of $\Gamma$ on the two sides.

Let us suppose moreover the generic fibre $X_\eta$ is smooth.
In order to compute the local factor
of the Hasse-Weil zeta function, one should calculate the trace
$$\sum_j (-1)^j{\mathrm Tr}({\mathrm Fr}_q,{\mathrm H}^j(X_{\bar\eta},{{\bar{\Bbb Q}}_\ell})^{\Gamma_0}).$$
Assuming that the graded pieces in the monodromy filtration of 
${\mathrm H}^j(X_{\bar\eta},{{\bar{\Bbb Q}}_\ell})$ are pure (Deligne's conjecture),
Rapoport proved that the true local factor 
is completely determined by the semi-simple local factor, see \cite{Rapoport}.
Now by the above discussion the semi-simple trace can be computed by the formula
$$\sum_j (-1)^j{\mathrm Tr}^{ss}({\mathrm Fr}_q,{\mathrm H}^j(X_{\bar\eta},{{\bar{\Bbb Q}}_\ell}))
=\sum_{x\in X(k)}{\mathrm Tr}^{ss}({\mathrm Fr}_q,{\mathrm R}\Psi({{\bar{\Bbb Q}}_\ell})_x).$$

\section{Statement of the main result}
\subsection{Nearby cycles on local models}

We have seen in subsection 2.3 (resp. 2.5 for symplectic case)  that the generic fibre of $M_{n_\pm}$
admits a stratification with smooth strata
$$ M_{n_\pm,\eta}=\coprod_{\lambda\in\Lambda(n_\pm)} O_{\lambda,\eta}.$$
Denote by ${\bar O}_{\lambda,\eta}$ the Zariski closure of $O_{\lambda,\eta}$ 
in $M_{n_\pm,\eta}$; in general ${\bar O}_{\lambda, \eta}$ is no longer smooth . 
It is natural to consider ${\mathcal A}_{\lambda,\eta}={\mathrm IC}(O_{\lambda,\eta})$, its $\ell$-adic intersection complex.

We want to calculate the function 
$$\tau^{ss}_{{\mathrm R}\Psi^M({\mathcal A}_{\lambda,\eta})}(x)
={\mathrm Tr}^{ss}({\mathrm Fr}_q,{\mathrm R}\Psi^M({\mathcal A}_{\lambda,\eta})_x)$$
of semi-simple trace of the Frobenius endomorphism on the nearby
cycle complex ${\mathrm R}\Psi^M({\mathcal A}_{\lambda,\eta})$ defined in the last section.  We are denoting the scheme $M_{n_\pm}$ simply by $M$ here.

As $O_{\lambda,\eta}$ is an orbit of $J_{n_\pm,\eta}$,
the intersection complex ${\mathcal A}_{\lambda,\eta}$ is naturally 
$J_{n_\pm,\eta}$-equivariant. As we know that $J_{n_\pm}$
is smooth over $S$ by Lemma 3, its nearby cycle
${\mathrm R}\Psi^M({\mathcal A}_{\lambda,\eta})$ is $J_{n_\pm,{\bar s}}$-equivariant.
In particular, the function
$$\tau^{ss}_{{\mathrm R}\Psi^M({\mathcal A}_{\lambda,\eta})}:M_{n_\pm}(k)\rightarrow{{\bar{\Bbb Q}}_\ell}$$
is $J_{n_\pm}(k)$-invariant.

Now following the group theoretic description of the action of $J_{n_\pm}(k)$
on $M_{n\pm}(k)$ in subsection 2.4 (resp. 2.5), 
we can consider the function $\tau^{ss}_{{\mathrm R}\Psi^M({\mathcal A}_{\lambda,\eta})}$
as a function on $G_k$ with compact support which is invariant on the left 
and on the right by the Iwahori subgroup $I_k$
$$\tau^{ss}_{{\mathrm R}\Psi^M({\mathcal A}_{\lambda,\eta})}\in{\mathcal H}(G_k/\!/I_k).$$

The following statement was conjectured by R. Kottwitz, and is the main result of this paper.

\begin{theorem}Let $G$ be either ${\mathrm GL}(d)$ or ${\mathrm GSp}(2d)$.  Let $M = M_{n_\pm}$ be the scheme associated to the group $G$ and the pair of integers $n_\pm$, as above.  Then we have the formula
$$\tau^{ss}_{{\mathrm R}\Psi^M({\mathcal A}_{\lambda,\eta})}=(-1)^{2\langle\rho,\lambda\rangle }
\sum_{\lambda'\leq\lambda}m_\lambda(\lambda') z_{\lambda'}$$
where $z_{\lambda'}$ is the function of Bernstein associated to the 
dominant coweight $\lambda'$, which lies in the center $Z({\mathcal H}(G_k/\!/I_k))$
of ${\mathcal H}(G_k/\!/I_k)$.
\end{theorem}

Here, $\rho$ is half the sum of positive roots for $G$ and thence $2\langle\rho,\lambda\rangle $
is the dimension of $O_{\lambda,\eta}$. The integer $m_\lambda(\lambda')$ 
is the multiplicity of weight $\lambda'$ occuring in the representation of
highest weight $\lambda$.  The partial ordering $\lambda' \leq \lambda$ is defined to mean that $\lambda - \lambda'$ is a sum of positive coroots of $G$.

Comparing with the formula for minuscule $\mu$ given in Kottwitz' conjecture (cf. Introduction), one notices
the absence of the factor $q^{\langle\rho,\mu\rangle }$ and the appearance
of the sign $(-1)^{2\langle\rho,\mu\rangle }$. This difference is explained by
the normalization of the intersection complex ${\mathcal A}_{\mu,\eta}$.
For minuscule coweights $\mu$, the orbit $O_\mu$ is closed. 
Consequently, the intersection complex ${\mathcal A}_{\mu,\eta}$ differs from the 
constant sheaf only by normalization factor
$${\mathcal A}_{\mu,\eta}={{\bar{\Bbb Q}}_\ell}[2\langle\rho,\mu\rangle ](\langle\rho,\mu\rangle ).$$

We refer to Lusztig's article \cite{Lusztig} for the definition of 
Bernstein's functions. In fact, what we need is rather the properties
that characterize these functions. We will recall these 
properties in the next subsection. 

\subsection{A commutative triangle}

Denote by $K_k$ the standard maximal compact subgroup $G(k[[t]])$ of $G_k$, where $G$ is either ${\mathrm GL}(d)$ or ${\mathrm GSp}(2d)$.
The ${{\bar{\Bbb Q}}_\ell}$-valued functions with compact support in $G_k$ invariant 
on the left and on the right by $K_k$ form a commutative algebra 
${\mathcal H}(G_k/\!/K_k)$ with respect to the convolution product. Here the convolution is defined using the Haar measure on $G_k$ which gives $K_k$ measure 1.  Denote 
by ${\Bbb I}_K$ the characteristic function of $K_k$. This element
is the unit of the algebra ${\mathcal H}(G_k/\!/K_k)$.  Similarly we define the convolution on ${\mathcal H}(G_k /\!/I_k)$ using the Haar measure on $G_k$ which gives $I_k$ measure 1.

We consider the following triangle
$$\diagramme{ &{{\bar{\Bbb Q}}_\ell}[X_*]^{W}  & \cr
\hfill{}^{\mathrm Bern.}\! \swarrow & & 
\nwarrow\!{}^{\mathrm Sat.}\hfill \cr
Z({\mathcal H}(G_k/\!/I_k)) & \hfld{}{-*{\Bbb I}_K} & {\mathcal H}(G_k/\!/K_k)}$$ 
Here ${{\bar{\Bbb Q}}_\ell}[X_*]^{W}$ is the $W$-invariant sub-algebra of the ${{\bar{\Bbb Q}}_\ell}$-algebra
associated to the group of cocharacters of the standard (diagonal) torus $T$ in $G$ and $W$ is the Weyl group associated to $T$. 
For the case $G = {\mathrm GL}(d)$, this algebra is isomorphic to the algebra of symmetric polynomials 
with $d$ variables and their inverses: ${\bar {\Bbb Q}}_{\ell}[X^{\pm}_1,\ldots,X^{\pm}_d]^{S_d}$. 

The above maps 
$${\mathrm Sat}:{\mathcal H}(G_k/\!/K_k)\rightarrow{{\bar{\Bbb Q}}_\ell}[X_*]^{W}$$
and 
$${\mathrm Bern}:{{\bar{\Bbb Q}}_\ell}[X_*]^{W}\rightarrow Z({\mathcal H}(G_k/\!/I_k))$$
are the isomorphisms of algebras constructed by Satake, see \cite{Satake}
and by Bernstein, see \cite{Lusztig}.  
It follows immediately from its definition that 
the Bernstein isomorphism sends the irreducible
character $\chi_\lambda$ of highest weight $\lambda$ to 
$${\mathrm Bern}(\chi_\lambda)=\sum_{\lambda'\leq\lambda}m_\lambda(\lambda')
z_{\lambda'}.$$ 

The horizontal map
$$Z({\mathcal H}(G_k/\!/I_k))\rightarrow {\mathcal H}(G_k/\!/K_k)$$
is defined by $f\mapsto f*{\Bbb I}_K$ where
$$f*{\Bbb I}_K(g)=\int_{G_k}f(gh^{-1}){\Bbb I}_K(h) \,d h.$$

The next statement seems to be known to the experts. It can be deduced 
easily, see \cite{Haines}, from results of Lusztig \cite{Lusztig} and 
Kato \cite{Kato}. Another proof can be found in an article of Dat \cite{Dat}.

\begin{lemma}
The above triangle is commutative.
\end{lemma}
It follows that the horizontal map is an isomorphism, and that \linebreak $(-1)^{2\langle \rho\, , \, \lambda \rangle} \sum_{\lambda' \leq \lambda}m_{\lambda}(\lambda')z_{\lambda'}$ is the unique element in $Z({\mathcal H}(G_k/\!/I_k))$ whose image in ${\mathcal H}(G_k/\!/K_k)$ has Satake transform $(-1)^{2\langle \rho\, , \, \lambda \rangle} \chi_{\lambda}$.

Thus in order to prove the Theorem 11, it suffices now to prove the two
following statements.

\begin{proposition}
The function $\tau^{ss}_{{\mathrm R}\Psi^M({\mathcal A}_{\lambda,\eta})}$ lies in the 
center $Z({\mathcal H}(G_k/\!/I_k))$ of the algebra ${\mathcal H}(G_k/\!/I_k)$.
\end{proposition}

\begin{proposition}
The Satake transform of 
$\tau^{ss}_{{\mathrm R}\Psi^M({\mathcal A}_{\lambda,\eta})}*{\Bbb I}_K$
is equal to $(-1)^{2\langle\rho,\lambda\rangle }\chi_\lambda$, where $\chi_\lambda$
is the irreducible character of highest weight $\lambda$.
\end{proposition}

In fact we can reformulate Proposition 14 in such a way that 
it becomes independent of Proposition  13. 
We will prove Proposition 14 in the next section. 

In order to prove Proposition 13, we have to adapt Lusztig's
construction of geometric convolution to our context. 
This will be done in the section 7. 
The proof of Proposition 14 itself will be given in section 8.

\section{Proof of Proposition 14}
\subsection{Averaging by $K$}

The map $$Z({\mathcal H}(G_k/\!/I_k))\rightarrow {\mathcal H}(G_k/\!/K_k)$$
defined by $f\mapsto f*{\Bbb I}_K$ can be obviously extended
to a map
$$C_c(G_k/I_k)\rightarrow C_c(G_k/K_k)$$
where $C_c(G_k/I_k)$ (resp. $C_c(G_k/K_k)$) is the space of functions
with compact support in $G_k$ invariant on the right by $I_k$ (resp. $K_k$).
This map can be rewritten as follows
$$f*{\Bbb I}_K(g)=\sum_{h\in K_k/I_k} f(gh).$$

Therefore, this operation corresponds to summing along the fibres of 
the map $G_k/I_k\rightarrow G_k/K_k$.
For the particular function $\tau^{ss}_{{\mathrm R}\Psi^M({\mathcal A}_{\lambda,\eta})}$, 
it amounts to summing along the fibres of the map
$$\pi(k):M_{n_\pm}(k)\rightarrow {\mathrm Grass}_{n_\pm}(k),$$
(see Lemma 6).

By using now the sheaf-function dictionary for semi-simple trace,
we get
$$\tau^{ss}_{{\mathrm R}\Psi^M({\mathcal A}_{\lambda,\eta})}*{\Bbb I}_K
=\tau^{ss}_{{\mathrm R}\pi_{{\bar s},*}{\mathrm R}\Psi^M({\mathcal A}_{\lambda,\eta})}.$$
The nearby cycle functor commutes with direct image by a proper morphism, 
so that 
$${\mathrm R}\pi_{{\bar s},*}{\mathrm R}\Psi^M({\mathcal A}_{\lambda,\eta})
={\mathrm R}\Psi^{{\mathrm Grass}}{\mathrm R}\pi_{\eta,*}({\mathcal A}_{\lambda,\eta}).$$
By Lemma 4, $\pi_{\eta}$ is an isomorphism. Consequently,
${\mathrm R}\pi_{\eta,*}({\mathcal A}_{\lambda,\eta})={\mathcal A}_{\lambda,\eta}$.

According to the description of ${\mathrm Grass} = {\mathrm Grass}_{n_\pm}$ (see subsections 2.3 and 2.5),
we can prove that ${\mathrm R}\Psi^{{\mathrm Grass}}{\mathcal A}_{\lambda,\eta}={\mathcal A}_{\lambda,{\bar s}}$ (note that the complex ${\mathcal A}_{\lambda,\eta}$ over ${\mathrm Grass}_{\eta}$ can be extended in a canonical fashion to a complex ${\mathcal A}_{\lambda}$ over the $S$-scheme ${\mathrm Grass}$, thus ${\mathcal A}_{\lambda, {\bar s}}$ makes sense).
In particular, the inertia subgroup $\Gamma_0$ acts trivially on 
${\mathrm R}\Psi^{{\mathrm Grass}}{\mathcal A}_{\lambda,\eta}$ and the semi-simple trace is just
the ordinary trace. The proof of a more general statement will be given in the 
following appendix. 

By putting together the above equalities, we obtain
$${\mathrm R}\pi_{{\bar s},*}{\mathrm R}\Psi^M({\mathcal A}_{\lambda,\eta})={\mathcal A}_{\lambda,s}.$$

To conclude the proof of Proposition 14, we quote an important
theorem of Lusztig and Kato, see \cite{Lusztig} and \cite{Kato}.
We remark that Ginzburg and also Mirkovic and Vilonen have put this result
in its natural framework : a Tannakian equivalence, 
see \cite{Ginzburg},\cite{Mirkovic-Vilonen}.

\begin{theorem} [Lusztig, Kato]
The Satake transform of the function $\tau^{ss}_{{\mathcal A}_{\lambda,s}}$ is equal to
$${\mathrm Sat}(\tau_{{\mathcal A}_{\lambda,s}})=(-1)^{2\langle\rho,\lambda\rangle }\chi_\lambda$$
where $\chi_\lambda$ is the irreducible character of highest weight 
$\lambda$. 
\end{theorem}

\subsection{Appendix}

This appendix seems to be well known to the experts.
We thank G. Laumon who has kindly explained it to us.

Let us consider the following situation.

Let $X$ be a proper scheme over $S$ equipped with 
an action of a group scheme $J$ smooth over $S$. 
We suppose there is a stratification
$$X=\coprod_{\alpha\in\Delta} X_\alpha$$
with each stratum $X_\alpha$ smooth over $S$. We assume that the group scheme $J$ acts 
transitively on all fibers of $X_\alpha$. Moreover, we suppose 
there exists, for  each $\alpha$, a $J$-equivariant resolution of singularities
${\tilde X}_\alpha$ 
$$\pi_\alpha:{\tilde X}_\alpha\rightarrow{\bar X}_\alpha$$
of the closure ${\bar X}_\alpha$ of $X_\alpha$, such 
that this resolution ${\tilde X}_\alpha$, smooth over $S$,
contains $X_\alpha$ as a Zariski open ;
the complement ${\tilde X}_\alpha-X_\alpha$ is also supposed to be a union
of normal crossing divisors. 

If $X$ is an invariant subscheme of the affine Grassmannian 
or of the affine flag variety, we can use the Demazure resolution.

Let $i_\alpha$ denote the inclusion map $X_\alpha\rightarrow X$
and let ${\mathcal F}_\alpha$ denote $i_{\alpha,!}{{\bar{\Bbb Q}}_\ell}$. 
A complex of sheaves 
${\mathcal F}$ is said $\Delta$-constant if
its cohomology sheaves of ${\mathcal F}$ are 
successive extensions of ${\mathcal F}_\alpha$ with $\alpha\in\Delta$. 
The intersection complex
of ${\bar X}_\alpha$ is $\Delta$-constant.

For an $\ell$-adic complex ${\mathcal F}$ of sheaves on $X$, 
there exists a canonical morphism
$${\mathcal F}_{\bar s}\rightarrow{\mathrm R}\Psi^X({\mathcal F}_\eta)$$ 
whose the mapping cylinder is the vanishing cycle
${\mathrm R}\Phi^X({\mathcal F})$. 

\begin{lemma} 
If ${\mathcal F}$ is $\Delta$-constant bounded complex, 
${\mathrm R}\Phi^X({\mathcal F})=0$.
\end{lemma}

\PROOF
Clearly, it suffices to prove ${\mathrm R}\Phi^X({\mathcal F}_\alpha)=0$.
Consider the equivariant resolution 
$\pi_\alpha:{\tilde X}_\alpha\rightarrow {\bar X}_\alpha$. 
We have a canonical isomorphism 
$${\mathrm R}\pi_{\alpha,*}{\mathrm R}\Phi^{{\tilde X}_\alpha}({\mathcal F}_\alpha)
\ident{\mathrm R}\Phi^{{\bar X}_\alpha}({\mathcal F}_\alpha).$$
It suffices then to prove 
${\mathrm R}\Phi^{{\tilde X}_\alpha}({\mathcal F}_\alpha)=0$.
This is known because ${\tilde X}_\alpha$ is smooth over $S$ and 
${{\tilde X}_\alpha}-X_\alpha$ is union of normal crossing divisors. $\square$

\begin{corollary}
If ${\mathcal F}$ is $\Delta$-constant and bounded,
the inertia group $\Gamma_0$ acts trivially on 
the nearby cycle ${\mathrm R}\Psi^X({\mathcal F}_\eta)$.
\end{corollary}

\PROOF
The morphism ${\mathcal F}_{\bar s}\rightarrow{\mathrm R}\Psi^X({\mathcal F}_\eta)$
is an isomorphism compatible with the actions of $\Gamma$. The
inertia subgroup $\Gamma_0$ acts trivially on ${\mathcal F}_{\bar s}$, thus 
it acts trivially on ${\mathrm R}\Psi^X({\mathcal F}_\eta)$, too. $\square$

\section{Invariant subschemes of $G/I$}

We recall here the well known ind-scheme structure of $G_k/I_k$
where $G$ denotes the group ${\mathrm GL}(d,k(\!(t+\varpi)\!))$ or 
the group ${\mathrm GSp}(2d,k(\!(t+\varpi)\!))$
and where $I$ is its standard Iwahori subgroup. 
The variable  $t+\varpi$ is used instead of $t$ in order to be compatible with
the definitions of local models given in section 2.

\subsection{Linear case}

Let $N_{n_\pm}$ be the functor which associates to each ${\mathcal O}$-algebra $R$
the set of 
$${\mathcal L}_\bullet=({\mathcal L}_0\subset{\mathcal L}_1\subset\cdots\subset{\mathcal L}_d=(t+\varpi)^{-1}{\mathcal L}_0)$$
where ${\mathcal L}_0,{\mathcal L}_1,\ldots$ are $R[t]$-submodules of $R[t,t^{-1},(t+\varpi)^{-1}]^d$
such that for $i=0,1,\ldots,d-1$
$$(t+\varpi)^{n_+}{\mathcal V}_{i,R}\subset {\mathcal L}_i\subset (t+\varpi)^{n_-}{\mathcal V}_{i,R}$$
and ${\mathcal L}_i / (t+\varpi)^{n_+}{\mathcal V}_{i,R}$ is locally a direct factor, of fixed rank independent of $i$, of the free $R$-module $(t+\varpi)^{n_-}{\mathcal V}_{i,R} / (t+\varpi)^{n_+}{\mathcal V}_{i,R}$. Obviously, this functor is represented by a closed subscheme in a product of Grassmannians. In particular, $N_{n_\pm}$ is proper.

Let $I_{n_\pm}$ be the functor which associates to each ${\mathcal O}$-algebra $R$
the group $R[t]$-linear automorphisms of 
$$  (t+\varpi)^{n_--1}R[t]^d/(t+\varpi)^{n_+}R[t]^d$$
fixing the image in this quotient of the filtration
$${\mathcal V}_{0,R}\subset{\mathcal V}_{1,R}\subset\cdots\subset {\mathcal V}_{d,R}=(t+\varpi)^{-1}{\mathcal V}_{0,R}.$$
This functor is represented by a smooth group scheme over $S$ which
acts on $N_{n_\pm}$.

\subsection{Symplectic case}

Let $N_{n_\pm}$ be the functor which associates to each 
${\mathcal O}$-algebra $R$ the set of sequences
$${\mathcal L}_\bullet=({\mathcal L}_0\subset{\mathcal L}_1\subset\cdots\subset{\mathcal L}_d)$$
where ${\mathcal L}_0,{\mathcal L}_1,\ldots$ are $R[t]$-submodules of 
$R[t,t^{-1},(t+\varpi)^{-1}]^{2d}$ satisfying
$$(t+\varpi)^{n_+}{\mathcal V}_{i,R}\subset {\mathcal L}_i\subset (t+\varpi)^{n_-}{\mathcal V}_{i,R}$$
and such that ${\mathcal L}_i / (t+\varpi)^{n_+}{\mathcal V}_{i,R}$ is locally a direct factor of $(t+\varpi)^{n_-}{\mathcal V}_{i,R} / (t+\varpi)^{n_+}{\mathcal V}_{i,R}$ of rank $(n_+ - n_-)d$ for all $i=0,1,\ldots,d$, and ${\mathcal L}_0$ (resp. ${\mathcal L}_d$) 
is autodual 
with respect to the symplectic form $(t+\varpi)^{-n_--n_+}\langle\,,\,\rangle $ 
(resp. $(t+\varpi)^{-n_--n_++1}\langle\,,\,\rangle $).

Let $I_{n_\pm}$ be the functor which associates to each ${\mathcal O}$-algebra $R$
the group $R[t]$-linear automorphisms of 
$$  (t+\varpi)^{n_--1}R[t]^{2d}/(t+\varpi)^{n_+}R[t]^{2d}$$
fixing the image in this quotient, of the filtration
$${\mathcal V}_{0,R}\subset{\mathcal V}_{1,R}\subset\cdots\subset {\mathcal V}_{2d,R}=(t+\varpi)^{-1}{\mathcal V}_{0,R},$$
and fixing the symplectic form $(t+\varpi)^{-n_--n_++1}\langle\,,\,\rangle $ up to a unit in $R$.
This functor is represented by a smooth group scheme over $S$ which
acts on $N_{n_\pm}$.

\subsection{There is no vanishing cycle on $N$}

It is well known (see \cite{Mathieu} for instance) that $N = N_{n_\pm}$ admits 
a stratification by $I_{n_\pm}$-orbits
$$N_{n_\pm}=\coprod_{w\in {\tilde W}'(n_\pm)} O_w$$
where ${\tilde W}'(n_\pm)$ is a finite subset of the affine Weyl group 
${\tilde W}$ of ${\mathrm GL}(d)$ (resp. ${\mathrm GSp}(2d)$). 
For all $w\in {\tilde W}'(n_\pm)$, $O_w$
is smooth over $S$  and $I_{n_\pm}$ acts transitively on its geometric fibers.
All this remains true if we replace $S$ by any other base scheme.

Let ${\bar O}_w$ denote the closure of $O_w$. Let ${\mathcal I}_{w,\eta}$
(resp. ${\mathcal I}_{w,s}$) denote the intersection complex of ${\bar O}_{w,\eta}$
(resp. ${\bar O}_{w,s}$). We have
$${\mathrm R}\Psi^{N}({\mathcal I}_{w,\eta})={\mathcal I}_{w,{\bar s}}$$
(see Appendix 5.2 for a proof). In particular, the inertia subgroup 
$\Gamma_0$ acts trivially on ${\mathrm R}\Psi^{N}({\mathcal I}_{w,\eta})$.

Let ${\tilde W}$ be the affine Weyl group of ${\mathrm GL}(d)$, respectively
${\mathrm GSp}(2d)$. It can be easily checked that
${\tilde W}=\bigcup_{n_\pm}\  {\tilde W}'(n_\pm)$
for the linear case as well as for the symplectic case. 

\section{Convolution product of ${\mathcal A}_\lambda$ with ${\mathcal I}_w$} 
\subsection{Convolution diagram}

In this section, we will adapt a
construction due to Lusztig in order 
to define the convolution product of an equivariant perverse sheaf 
${\mathcal A}_\lambda$ over $M_{n_\pm}$ with an equivariant perverse 
sheaf ${\mathcal I}_w$ over $N_{n'_\pm}$.
See Lusztig's article \cite{Lusztig1} for a quite general construction.

For any dominant coweight $\lambda$ and any $w\in{\tilde W}$,
we can choose $n_\pm$ and $n'_\pm$ so that
$\lambda\in\Lambda(n_\pm)$ and $w\in{\tilde W}'(n'_\pm)$.  
>From now on, since $\lambda$ and $w$ as well as $n_\pm$ and $n'_\pm$ 
are fixed,
we will often write $M$ for $M_{n_\pm}$ and $N$ for $N_{n'_\pm}$.
This should not cause any confusion.

The aim of this subsection is to construct the convolution 
diagram \`a la Lusztig
$$\diagramme{ & {\tilde M}\times {\tilde N} & & &\cr
\hfill{}^{p_1}\! \swarrow & & 
\searrow\!{}^{p_2}\hfill & & \cr
M\times N & & M\,{\tilde\times}\,N
&\hfld{m}{} & P}$$
with the usual properties that will be made precise later.

\subsection{Linear case}

\begin{itemize}
\item
The functor $M\,{\tilde\times}\,N$ associates to each 
${\mathcal O}$-algebra $R$ the set of pairs $({\mathcal L}_\bullet,{\mathcal L}'_\bullet)$
$$\displaylines{
{\mathcal L}_\bullet=({\mathcal L}_0\subset{\mathcal L}_1\subset\cdots\subset{\mathcal L}_d=(t+\varpi)^{-1}{\mathcal L}_0)\cr
{\mathcal L}'_\bullet=({\mathcal L}'_0\subset{\mathcal L}'_1\subset\cdots\subset{\mathcal L}'_d=(t+\varpi)^{-1}{\mathcal L}'_0)
}$$
where ${\mathcal L}_i,{\mathcal L}'_i$ are $R[t]$-submodules of $R[t,t^{-1},(t+\varpi)^{-1}]^d$
satisfying the following conditions
$$\displaylines{
t^{n_+}{\mathcal V}_{i,R}\subset{\mathcal L}_i\subset t^{n_-}{\mathcal V}_{i,R}\cr
(t+\varpi)^{n'_+}{\mathcal L}_{i}\subset{\mathcal L}'_i\subset (t+\varpi)^{n'_-}{\mathcal L}_{i}}$$
As usual, ${\mathcal L}_i/t^{n_+}{\mathcal V}_{i,R}$ is supposed to be locally a direct factor of $t^{n_-}{\mathcal V}_{i,R}/t^{n_+}{\mathcal V}_{i,R}$, and ${\mathcal L}'_i/(t+\varpi)^{n'_+}{\mathcal L}_{i}$ locally a direct factor of $(t+\varpi)^{n'_-}{\mathcal L}_{i}/(t+\varpi)^{n'_+}{\mathcal L}_{i}$ as $R$-modules.The ranks of the projective $R$-modules ${\mathcal L}_i/t^{n_+}{\mathcal V}_{i,R}$ and ${\mathcal L}'_i/(t+\varpi)^{n'_+}{\mathcal L}_{i}$ are each also supposed to be independent of $i$. It follows from the above conditions that   
$$t^{n_+}(t+\varpi)^{n'_+}{\mathcal V}_{i,R}\subset{\mathcal L}'_i\subset t^{n_-}(t+\varpi)^{n'_-}{\mathcal V}_{i,R}$$
and ${\mathcal L}'_i/t^{n_+}(t+\varpi)^{n'_+}{\mathcal V}_{i,R}$ is locally a direct factor of $t^{n_-}(t+\varpi)^{n'_-}{\mathcal V}_{i,R}/t^{n_+}(t+\varpi)^{n'_+}{\mathcal V}_{i,R}$ as an $R$-module. Thus defined the functor $M \, {\tilde \times} \, N$ is represented by a projective scheme over $S$.

\item
The functor $P$ associates to each 
${\mathcal O}$-algebra $R$ the set of chains ${\mathcal L}'_\bullet$
$${\mathcal L}'_\bullet=({\mathcal L}'_0\subset{\mathcal L}'_1\subset\cdots\subset{\mathcal L}'_d=(t+\varpi)^{-1}{\mathcal L}'_0)
$$
where ${\mathcal L}'_i$ are $R[t]$-submodules of $R[t,t^{-1},(t+\varpi)^{-1}]^d$
satisfying 
$$t^{n_+}(t+\varpi)^{n'_+}{\mathcal V}_{i,R}\subset{\mathcal L}'_i\subset 
t^{n_-}(t+\varpi)^{n'_-}{\mathcal V}_{i,R}$$
and the usual conditions ``locally a direct factor as $R$-modules''.
As above, ${\rm rk}_R({\mathcal L}'_i/t^{n_+}(t+\varpi)^{n'_+}{\mathcal V}_{i,R})$ 
is supposed to be independent of $i$. Obviously, this functor is represented by a projective scheme over $S$.

\item
The forgetting map $m({\mathcal L}_\bullet,{\mathcal L}'_\bullet)={\mathcal L}'_\bullet$
yields a morphism
$$m:M\,{\tilde\times}\,N\rightarrow P.$$
This map is defined: it suffices to note that 
$t^{n_-}(t+\varpi)^{n'_-}{\mathcal V}_{i,R} / {\mathcal L}'_i$ is locally free as an $R$-module, 
being an extension of $t^{n_-}{\mathcal V}_{i,R} / {\mathcal L}_i$ by 
$(t+\varpi)^{n'_-}{\mathcal L}_i / {\mathcal L}'_i$, each of which is locally free. 
Clearly, this morphism is a proper morphism because its source and its target 
are proper schemes over $S$. 
\end{itemize}

Now before we can construct the schemes ${\tilde M}$, ${\tilde N}$, and the remaining morphisms in the convolution diagram, we need the following simple remark.

\begin{lemma}
The functor which associates to each ${\cal O}$-algebra $R$ the set of matrices $g\in{\frak{gl}}_s(R)$ such that the image of $g:R^s\rightarrow R^s$ is locally a direct factor of rank $r$ of $R^s$ is representable by a locally closed subscheme of ${\frak{gl}}_s$.
\end{lemma}

\PROOF
For $1\leq i\leq s$, denote by ${\rm St}_i$ the closed subscheme of ${\frak{gl}}_s$ defined by the equations : all minors of order at least $i+1$ vanish. 
By using Nakayama's lemma, one can see easily that the above functor is 
represented by the quasi-affine, locally closed subscheme ${\rm St}_r-{\rm St}_{r-1}$ of ${\frak{gl}}_s$. $\square$

\bigskip
Now let ${\bar{\mathcal V}}_0\subset{\bar{\mathcal V}}_1\subset\cdots$ be the image of
${\mathcal V}_0\subset{\mathcal V}_1\subset\cdots$ in the quotient 
$${\bar{\mathcal V}}=t^{n_-}(t+\varpi)^{n'_--1}{\mathcal O}[t]^d/
t^{n_+}(t+\varpi)^{n'_+}{\mathcal O}[t]^d.$$
Let ${\bar{\mathcal L}}_0\subset{\bar{\mathcal L}}_1\subset\cdots$ be the images of 
${\mathcal L}_0\subset{\mathcal L}_1\subset\cdots$ in the quotient
${\bar{\mathcal V}}_{R}={\bar{\mathcal V}}\otimes_{\mathcal O} R$.
Because ${\mathcal L}_i$ is completely determined by ${\bar{\mathcal L}}_i$, we can write
${\bar{\mathcal L}}_\bullet\in M(R)$ for ${\mathcal L}_\bullet\in M(R)$ and so on.

\begin{itemize}
\item
We consider the functor ${\tilde M}$
which associates to each ${\mathcal O}$-algebra $R$
the set of $R[t]$-endomorphisms
$g\in{\mathrm End}({\bar{\mathcal V}}_{R})$
such that if 
${\bar{\mathcal L}}_i=g\bigl(t^{n_-}{\bar{\mathcal V}}_i\bigr)$
then 
$$t^{n_+}{\bar{\mathcal V}}_{i,R}\subset{\bar{\mathcal L}}_i\subset
t^{n_-}{\bar{\mathcal V}}_{i,R}$$
and ${\bar{\mathcal L}}_i/t^{n_+}{\bar{\mathcal V}}_{i,R}$ is locally a direct factor
of $t^{n_-}{\bar{\mathcal V}}_{i,R}/t^{n_+}{\bar{\mathcal V}}_{i,R}$, of the same rank,
for all $i=0,\ldots,d-1$. 
Using Lemma 18 ones sees this functor is representable and comes naturally 
with a morphism $p:{\tilde M}\rightarrow M$.

\item
In a totally analogous way, we consider the functor ${\tilde N}$
which associates to each ${\mathcal O}$-algebra $R$
the set of $R[t]$-endomorphisms $g\in{\mathrm End}({\bar{\mathcal V}}_R)$
such that if 
${\bar{\mathcal L}}_i=g\bigl((t+\varpi)^{n'_-}{\bar{\mathcal V}}_{i,R}\bigr)$ 
then 
$$(t+\varpi)^{n'_+}{\bar{\mathcal V}}_{i,R}\subset{\bar{\mathcal L}}_i\subset
(t+\varpi)^{n'_-}{\bar{\mathcal V}}_{i,R}$$
and ${\bar{\mathcal L}}_i/(t+\varpi)^{n'_+}{\bar{\mathcal V}}_{i,R}$
is locally a direct factor of $(t+\varpi)^{n'_-}{\bar{\mathcal V}}_{i,R}/(t+\varpi)^{n'_+}{\bar{\mathcal V}}_{i,R}$, of the same rank for all $i=0,\ldots,d-1$. As above, the representability follows from Lemma 18. This functor comes naturally with a morphism $p':{\tilde N}\rightarrow N$.

\item
Now we define the morphism
$p_1:{\tilde M}\times{\tilde N}
\rightarrow {M}\times N$
by $p_1=p\times p'$.

\item
We define the morphism
$p_2:{\tilde M}\times{\tilde N}
\rightarrow {M}\,{\tilde\times}\, N$
by $p_2(g,g')=({\mathcal L}_\bullet,{\mathcal L}'_\bullet)$ with
$$({\mathcal L}_\bullet,{\mathcal L}'_\bullet)=(g(t^{n_-}{\mathcal V}_\bullet),
gg'(t^{n_-}(t+\varpi)^{n'_-}{\mathcal V}_\bullet)).$$
\end{itemize}

We have now achieved the construction of the convolution 
diagram. We need to prove some usual facts related to this diagram.  

\begin{lemma}
The morphisms $p_1$ and $p_2$ are smooth and surjective. 
Their relative dimensions are equal.
\end{lemma}

\PROOF
The proof is very similar to that of Lemma 3. 
Let us note that the morphism
$p:{\tilde M}\rightarrow M$ can be factored as $p=f\circ j$
where $j:{\tilde M}\rightarrow U$ is an open immersion 
and $f:U\rightarrow M$ is the vector bundle defined as follows.
For any ${\mathcal O}$-algebra $R$ and any ${\mathcal L}_\bullet\in M(R)$,
the fibre of $U$ over ${\mathcal L}_\bullet$ is the $R$-module
$$U({\mathcal L}_\bullet)=\bigoplus_{i=0}^{d-1}
(t+\varpi)^{n'_-}{\mathcal L}_i/t^{n_+}(t+\varpi)^{n'_+}{\mathcal V}_{i,R}.$$ 
The morphisms $p',p_1$ and $p_2$ can be described 
in the same manner. The equality of relative dimensions of $p_1$ and $p_2$ follows from Lemma 23 (proved in section 8) and the fact that they are each smooth. $\square$

\bigskip
Just as in subsection 2.2, we can consider the group valued 
functor ${\tilde J}$ which associates to each ${\mathcal O}$-algebra $R$ 
the group of $R[t]$-linear automorphisms of ${\bar{\mathcal V}}_{R}$
which fix the filtration 
${\bar{\mathcal V}}_0\subset{\bar{\mathcal V}}_1\subset\cdots\subset{\bar{\mathcal V}}_d$.
Obviously, this functor is represented by an affine algebraic group
scheme over $S$. The same proof as that of Lemma 3 proves 
that ${\tilde J}$ is smooth over $S$.  
Moreover, there are canonical morphisms of $S$-group schemes 
${\tilde J} \rightarrow J$ and ${\tilde J} \rightarrow I$, where $J = J_{n_\pm}$ (resp. $I = I_{n'_\pm}$) is the group scheme defined in subsection 2.2 (resp. 6.1).    

\begin{itemize}
\item
We consider the action $\alpha_1$ of 
${\tilde J} \times {\tilde J}$
on ${\tilde M}\times{\tilde N}$ defined by
$$\alpha_1(h,h';g,g')=(gh^{-1},g'{h'}^{-1}).$$
Clearly, this action leaves stable the fibres of 
$p_1:{\tilde M}\times{\tilde N}\rightarrow
{M}\times{N}$.

\item
We also consider the action $\alpha_2$ of 
${\tilde J}\times {\tilde J}$
on the same ${\tilde M}\times{\tilde N}$ defined by
$$\alpha_1(h,h';g,g')=(gh^{-1},hg'{h'}^{-1}).$$
Clearly, this action leaves stable the fibres of 
$p_2:{\tilde M}\times{\tilde N}\rightarrow
{M}\,{\tilde\times}\,{N}$.
\end{itemize}

\begin{lemma}
The action $\alpha_1$, respectively $\alpha_2$, is transitive on all 
geometric fibres of $p_1$, respectively $p_2$.  The geometric fibers of $p_1$, respectively $p_2$, are therefore connected.
\end{lemma}

\PROOF
Let $E$ be a (separably closed) field containing the fraction field $F$
of ${\mathcal O}$ or its residue field $k$. Let $g,g'$ be elements of 
${\tilde M}(E)$ such that 
$${\mathcal L}_\bullet=p(g)=p(g')\in M(E).$$

For all $i=0,\ldots,d-1$, denote by ${\hat{\mathcal V}}_i$ and ${\hat{\mathcal L}}_i$ the tensors
$$\displaylines{
{\hat{\mathcal V}}_i={\mathcal V}_i\otimes_{{\mathcal O}[t]} E[t]_{(t(t+\varpi))}\cr
{\hat{\mathcal L}}_i={\mathcal L}_i\otimes_{E[t]} E[t]_{(t(t+\varpi))}}
$$ 
where $E[t]_{(t(t+\varpi))}$ is the localized ring of $E[t]$
at the ideal $(t(t+\varpi))$, i.e., the ring ${\mathcal S}^{-1}E[t]$ where ${\mathcal S} = E[t] - \{(t) \cup (t+\varpi)\}$; this is a semi-local ring.
Of course, we can consider the modules ${\hat{\mathcal V}}_i$ and ${\hat{\mathcal L}}_i$
as $E[t]_{(t(t+\varpi))}$-submodules of $E(t)^d$.

Clearly, we have an isomorphism
$${\bar{\mathcal V}}_{E}=
t^{n_-}(t+\varpi)^{n'_--1}{\hat{\mathcal V}}_0/t^{n_+}(t+\varpi)^{n'_+}{\hat{\mathcal V}}_0$$
so that $E[t]$-endomorphisms of ${\bar{\mathcal V}}_{E}$
are the same as $E[t]_{(t(t+\varpi))}$-endomorphisms of
${\hat{\mathcal V}}_0$ taken modulo $t^{n_+-n_-}(t+\varpi)^{n'_+-n'_-+1}$.

By using the Nakayama lemma, $g$ and $g'$ can be lifted to 
$${\hat g}, {\hat g}'\in{\mathrm GL}(d,E(t))$$
such that 
$${\hat{\mathcal L}}_i={\hat g}t^{n_-}{\hat{\mathcal V}}_i\ ; \ 
{\hat{\mathcal L}}_i={\hat g}'t^{n_-}{\hat{\mathcal V}}_i.$$
This induces of course ${\hat h}{\bar{\mathcal V}}_i={\bar{\mathcal V}}_i$
with ${\hat h}={\hat g}^{-1}{\hat g}'$ and  for all $i=0,\ldots,d-1$.

Let $h$ be the reduction modulo $t^{n_+-n_-}(t+\varpi)^{n'_+-n'_-+1}$
of ${\hat h}$. It is clear that $g'=gh$ and $h$ lies in ${\tilde J}(E)$.

We have proved that ${\tilde J}$ acts transitively on the geometric 
fibres of ${\tilde M}\rightarrow M$. We can prove in 
a completely similar way that ${\tilde J}$ acts transitively 
on geo\-metric  fibres of ${\tilde N}\rightarrow N$.
Consequently, the action $\alpha_1$ is transitive on geometric fibres 
of $p_1$.

The proof of the statement for $\alpha_2$ and $p_2$ is similar. 
$\square$

\bigskip
The symmetric construction yields the following diagram
$$\diagramme{ & {\tilde N}\times {\tilde M} & & &\cr
\hfill{}^{p'_1}\! \swarrow & & 
\searrow\!{}^{p'_2}\hfill & & \cr
N\times M & & N\,{\tilde\times}\,M
&\hfld{m'}{} & P}$$
enjoying the same structures and properties.  More precisely, we define $N \, {\tilde \times} \, M$ as follows: for each ${\mathcal O}$-algebra $R$, let $(N \, {\tilde \times} \, M)(R)$ be the set of pairs $({\mathcal L}'_\bullet,{\mathcal L}_\bullet)$ 
$$\displaylines{
{\mathcal L}'_\bullet=({\mathcal L}'_0\subset{\mathcal L}'_1\subset\cdots\subset{\mathcal L}'_d=(t+\varpi)^{-1}{\mathcal L}'_0)\cr
{\mathcal L}_\bullet=({\mathcal L}_0\subset{\mathcal L}_1\subset\cdots\subset{\mathcal L}_d=(t+\varpi)^{-1}{\mathcal L}_0)
}$$
where ${\mathcal L}'_i,{\mathcal L}_i$ are $R[t]$-submodules of $R[t,t^{-1},(t+\varpi)^{-1}]^d$
satisfying the following conditions
$$\displaylines{
(t+\varpi)^{n'_+}{\mathcal V}_{i,R}\subset{\mathcal L}'_i\subset (t+\varpi)^{n'_-}{\mathcal V}_{i,R}\cr
t^{n_+}{\mathcal L}'_i\subset{\mathcal L}_i\subset t^{n_-}{\mathcal L}'_i}$$
such that for each $i=0,\ldots,d-1$, the $R$-module ${\mathcal L}'_i / (t+\varpi)^{n'_+}{\mathcal V}_{i,R}$ is locally a direct factor of $(t+\varpi)^{n'_-}{\mathcal V}_{i,R} / (t+\varpi)^{n'_+}{\mathcal V}_{i,R}$, and the $R$-module ${\mathcal L}_i/t^{n_+}{\mathcal L}'_i$ is locally a direct factor of $t^{n_-}{\mathcal L}'_i/t^{n_+}{\mathcal L}'_i$. It is also supposed that ${\rm rk}_R({\mathcal L}'_i / (t+\varpi)^{n'_+}{\mathcal V}_{i,R})$ and ${\rm rk}_R({\mathcal L}_i /t^{n_+}{\mathcal L}'_i)$ are independent of $i$.

The morphisms $p'_1$, $p'_2$, and $m'$ are defined in the obvious way: $p'_1 = p' \times p$, $m'({\mathcal L}'_\bullet, {\mathcal L}_\bullet) = {\mathcal L}_\bullet$, and $p'_2(g',g) = (g'(t+\varpi)^{n'_-}{\mathcal V}_{i,R} \, , \, g'g(t^{n_-}(t+\varpi)^{n'_-}){\mathcal V}_{i,R})$. 

\subsection{Symplectic case}

In this section we construct the symplectic analog of the convolution diagram just discussed.  In particular we need to define the schemes $M \, {\tilde \times} \, N$, ${\tilde M}$, 
${\tilde N}$, $P$, and the morphisms $p_1$, $p_2$, and $m$.  Moreover we need to construct the smooth group scheme ${\tilde J}$ which acts on the whole convolution diagram.  
Once this is done, defining the symplectic analogues of the actions $\alpha_1$ and $\alpha_2$, proving the symplectic analogues of Lemmas 19 and 20, and defining the symmetric construction are all straightforward tasks and will be left to the reader.

\begin{itemize}
\item
The functor $M\,{\tilde\times}\,N$ associates to each 
${\mathcal O}$-algebra $R$ the set of pairs $({\mathcal L}_\bullet,{\mathcal L}'_\bullet)$
$$\displaylines{
{\mathcal L}_\bullet=({\mathcal L}_0\subset{\mathcal L}_1\subset\cdots\subset{\mathcal L}_d)\cr
{\mathcal L}'_\bullet=({\mathcal L}'_0\subset{\mathcal L}'_1\subset\cdots\subset{\mathcal L}'_d)}$$
where ${\mathcal L}_i,{\mathcal L}'_i$ are $R[t]$-submodules of $R[t,t^{-1},(t+\varpi)^{-1}]^{2d}$
satisfying the following conditions
$$\displaylines{
t^{n_+}{\mathcal V}_{i,R}\subset{\mathcal L}_i\subset t^{n_-}{\mathcal V}_{i,R}\cr
(t+\varpi)^{n'_+}{\mathcal L}_{i}\subset{\mathcal L}'_i\subset (t+\varpi)^{n'_-}{\mathcal L}_{i}}$$
satisfying the usual ''locally direct factors as $R$-modules'' conditions: 
${\mathcal L}_i / t^{n_+}{\mathcal V}_{i,R}$ is locally a direct factor of $t^{n_-}{\mathcal V}_{i,R} / t^{n_+}{\mathcal V}_{i,R}$ of rank $(n_+ - n_-)d$ and ${\mathcal L}'_i / (t+\varpi)^{n'_+}{\mathcal L}_i$ is locally a direct factor of $(t+\varpi)^{n'_-}{\mathcal L}_i / (t+\varpi)^{n'_+}{\mathcal L}_i$ of rank $(n'_+ - n'_-)d$.
Moreover we suppose ${\mathcal L}_0$, ${\mathcal L}_d$, ${\mathcal L}'_0$ and ${\mathcal L}'_d$ are autodual with respect to
$t^{-n_--n_+}\langle\,,\,\rangle $, $t^{-n_--n_+}(t+\varpi)\langle\,,\,\rangle $,
$t^{-n_--n_+}(t+\varpi)^{-n'_--n'_+}\langle\,,\,\rangle $ and
$t^{-n_--n_+}(t+\varpi)^{-n'_--n'_++1}\langle\,,\,\rangle $ respectively.

\item
The functor $P$ associates to each 
${\mathcal O}$-algebra $R$ the set of chains ${\mathcal L}'_\bullet$
$${\mathcal L}'_\bullet=({\mathcal L}'_0\subset{\mathcal L}'_1\subset\cdots\subset{\mathcal L}'_d)$$
where ${\mathcal L}'_i$ are $R[t]$-submodules of $R[t,t^{-1},(t+\varpi)^{-1}]^{2d}$
satisfying 
$$t^{n_+}(t+\varpi)^{n'_+}{\mathcal V}_{i,R}\subset{\mathcal L}'_i\subset 
t^{n_-}(t+\varpi)^{n'_-}{\mathcal V}_{i,R},$$
such that the usual ``locally a direct factor as $R$-modules of rank $(n_+ - n_- + n'_+ - n'_-)d$'' condition holds,
and such that  ${\mathcal L}'_0$ and ${\mathcal L}'_d$ are autodual with respect to 
\linebreak $t^{-n_--n_+}(t+\varpi)^{-n'_--n'_+}\langle\,,\,\rangle $ and
$t^{-n_--n_+}(t+\varpi)^{-n'_--n'_++1}\langle\,,\,\rangle $ respectively.

\item
The forgetting map $m({\mathcal L}_\bullet,{\mathcal L}'_\bullet)={\mathcal L}'_\bullet$
yields a morphism
$m:M\,{\tilde\times}\,N\rightarrow P$.
Clearly, $m$ is a proper morphism between proper $S$-schemes. 

\item
We consider the functor ${\tilde M}$
which associates to each ${\mathcal O}$-algebra $R$
the set of $R[t]$-endomorphisms $g$ of
$${\bar{\mathcal V}}_{R}=t^{n_-}(t+\varpi)^{n'_--1}{\mathcal V}_{0,R}/
t^{n_+}(t+\varpi)^{n'_+}{\mathcal V}_{0,R}$$
satisfying
$$\langle gx,gy\rangle = c_g t^{n_+-n_-}\langle x,y\rangle $$
for some $c_g \in R^{\times}$, and such that if ${\bar{\mathcal L}}_i=g\bigl(t^{n_-}{\bar {\mathcal V}}_i\bigr)$ 
for $i=0,\ldots,d$, then we have 
$$t^{n_+}{\bar{\mathcal V}}_{i,R}\subset{\bar{\mathcal L}}_i\subset
t^{n_-}{\bar{\mathcal V}}_{i,R},$$
and ${\bar {\mathcal L}}_i / t^{n_+}{\bar {\mathcal V}}_{i,R}$ is locally a direct factor of $t^{n_-}{\bar {\mathcal V}}_{i,R} / t^{n_+}{\bar {\mathcal V}}_{i,R}$ of rank $(n_+ - n_-)d$. 
If $g\in {\tilde M}(R)$ then ones sees using the definitions that 
automatically,
${\bar{\mathcal L}}_\bullet=gt^{n_-}{\mathcal V}_{\bullet,R}\in M(R)$.
The functor ${\tilde M}$ is representable 
and comes naturally with a morphism $p:{\tilde M}\rightarrow M$.

\item
Next consider the functor ${\tilde N}$ which associates to each ${\mathcal O}$-algebra $R$ the set of $R[t]$-endomorphisms $g$ of ${\bar {\mathcal V}}_R$ satisfying
$$\langle gx,gy\rangle = c_g (t+\varpi)^{n'_+ - n'_-}\langle x,y \rangle$$
for some $c_g \in R^{\times}$ and such that
if ${\bar {\mathcal L}}'_i = g(t+\varpi)^{n'_-}{\bar {\mathcal V}}_{i,R}$ for $i=0,\ldots,d$ then we have
$$
(t+\varpi)^{n'_+}{\bar {\mathcal V}}_{i,R} \subset {\bar {\mathcal L}}'_i \subset (t+\varpi)^{n'_-}{\bar {\mathcal V}}_{i,R},
$$
and ${\bar {\mathcal L}}'_i / (t+\varpi)^{n'_+}{\bar {\mathcal V}}_{i,R}$ is locally a direct factor of $(t+\varpi)^{n'_-}{\bar {\mathcal V}}_{i,R} / (t+\varpi)^{n'_+}{\bar {\mathcal V}}_{i,R}$ of rank $(n'_+ - n'_-)d$.
>From the definitions one sees that ${\bar {\mathcal L}}'_\bullet \in N(R)$.  The functor ${\tilde N}$ is representable and comes with a functor $p': {\tilde N} \rightarrow N$.

\item We define $p_1 = p \times p'$.  We define $p_2 : {\tilde M} \times {\tilde N} \rightarrow M {\tilde \times} N$ exactly as in the linear case.   

\item We let ${\tilde J}$ denote the functor which associates to any ${\mathcal O}$-algebra $R$ the group of $R[t]$-linear automorphisms of ${\bar {\mathcal V}}_{R}$ which fix the form $t^{-n_- - n_+}(t+\varpi)^{-n'_- - n'_+ + 1}\langle \, , \, \rangle$ up to an element in $R^{\times}$.  As in Lemma 3, the group scheme ${\tilde J}$ is smooth over $S$.  There are canonical $S$-group scheme morphisms ${\tilde J} \rightarrow J$ and ${\tilde J} \rightarrow I$, where $J = J_{n_\pm}$ (resp. $I = I_{n'_\pm}$) was defined in subsection 2.5 (resp. 6.2).

\end{itemize}

\subsection{Definition of the convolution product}

Let us recall the standard definition of convolution product due to
Lusztig \cite{Lusztig1} (see also \cite{Ginzburg} and 
\cite{Mirkovic-Vilonen}).

Let $E$ be a field containing the fraction field $F$
of ${\mathcal O}$ or its residue field $k$ and let $\epsilon={\mathrm Spec\,}(E)\rightarrow S$
be the corresponding morphism. For all $S$-schemes $X$, let $X_\epsilon$
denote the base change $X\times_S\epsilon$.

Let ${\mathcal A}$ be a perverse sheaf over $M_{\epsilon}$ that is 
$J_{\epsilon}$-equivariant. Let ${\mathcal I}$ be a perverse sheaf 
over $N_{\epsilon}$ that is $I_{\epsilon}$-equivariant.
Both $I_{\epsilon}$ and $J_{\epsilon}$ are quotients of ${\tilde J}_{\epsilon}$, 
so we can say that 
${\mathcal A}$ and ${\mathcal I}$ are ${\tilde J}_{\epsilon}$-equivariant.

Since $p_1$ is a smooth morphism, the pull-back 
$p_1^*({\mathcal A}\boxtimes_\epsilon{\mathcal I})$ is also perverse up to the shift 
by the relative dimension of $p_1$. A priori, this pull-back is only 
$\alpha_1$-equivariant. As ${\mathcal A}$ and ${\mathcal I}$ are 
${\tilde J}_{\epsilon}$-equivariant, $p_1^*({\mathcal A}\boxtimes_\epsilon{\mathcal I})$
is also $\alpha_2$-equivariant.
Since $p_2$ is smooth and the action $\alpha_2$ is transitive on its 
geometric fibres, there exists a perverse sheaf 
${\mathcal A}\,{\tilde\boxtimes}_\epsilon{\mathcal I}$, unique up to unique isomorphism, 
such that
$$p_1^*({\mathcal A}\boxtimes_\epsilon{\mathcal I})=p_2^*({\mathcal A}\,{\tilde\boxtimes}_\epsilon{\mathcal I})$$
by the theorem 4.2.5 of Beilinson-Bernstein-Deligne \cite{BBD}.
And we put now
$${\mathcal A}*_\epsilon{\mathcal I}=
{\mathrm R} m_{*}({\mathcal A}\,{\tilde\boxtimes}_\epsilon{\mathcal I}).$$
By the symmetric construction, we can define the convolution product
${\mathcal I}*_\epsilon{\mathcal A}$.

Let $E$ be now the algebraic closure ${\bar k}$ of the residual field
$k$. We suppose that the perverse sheaves ${\mathcal A}$ and ${\mathcal I}$ are equipped with
an action of ${\mathrm Gal}({\bar F}/F)$ compatible with the action of 
${\mathrm Gal}({\bar F}/F)$ on the geometric special fibre through ${\mathrm Gal}({\bar k}/k)$.
In practice, the inertia subgroup $\Gamma_0$ acts trivially on ${\mathcal I}$ and 
non trivially on ${\mathcal A}$. As the semi-simple trace provides a 
sheaf-function dictionary, we have :
$$\displaylines{
\tau^{ss}_{\mathcal A}*\tau^{ss}_{\mathcal I}=\tau^{ss}_{{\mathcal A}\,*_{\bar s}\,{\mathcal I}}\cr
\tau^{ss}_{\mathcal I}*\tau^{ss}_{\mathcal A}=\tau^{ss}_{{\mathcal I}\,*_{\bar s}\,{\mathcal A}}}$$ 
where the convolution on the left hand is the ordinary convolution 
in the Hecke algebra ${\mathcal H}(G_k/\!/I_k)$.

\section{Proof of Proposition 13}

\subsection{Cohomological part}

According to the sheaf-function dictionary 
for semi-simple traces, it suffices to prove the following statement.
Beilinson and Gaitsgory have proved a related result in the equal characteristic case, using a deformation of the affine Grassmanian of $G$, see \cite{Gaitsgory}.

\begin{proposition}
We have an isomorphism
$${\mathrm R}\Psi^M({\mathcal A}_{\lambda,\eta})*_{\bar s}{\mathcal I}_{w,{\bar s}}\ident
{\mathcal I}_{w,{\bar s}}*_{\bar s}{\mathrm R}\Psi^M({\mathcal A}_{\lambda,\eta}).$$
\end{proposition}   

\PROOF
The above statement makes sense because the functor ${\mathrm R}\Psi$
sends perverse sheaves to perverse sheaves, by a theorem of Gabber,
see \cite{Illusie}.
In particular, ${\mathrm R}\Psi^M({\mathcal A}_{\lambda,\eta})$ is a perverse sheaf.

Let us recall that
$${\mathrm R}\Psi^N({\mathcal I}_{w,\eta})\ident{\mathcal I}_{w,s}$$
so that we have to prove
$${\mathrm R}\Psi^M({\mathcal A}_{\lambda,\eta})*_{\bar s}{\mathrm R}\Psi^N({\mathcal I}_{w,\eta})\ident
{\mathrm R}\Psi^N({\mathcal I}_{w,\eta})*_{\bar s}{\mathrm R}\Psi^M({\mathcal A}_{\lambda,\eta}).$$

First, let us prove that nearby cycle commutes with convolution
product. 

\begin{lemma}
We have the isomorphisms
$$\displaylines{
{\mathrm R}\Psi^M({\mathcal A}_{\lambda,\eta})*_{\bar s}{\mathrm R}\Psi^N({\mathcal I}_{w,\eta})\ident
{\mathrm R}\Psi^P({\mathcal A}_{\lambda,\eta}*_\eta{\mathcal I}_{w,\eta})\cr
{\mathrm R}\Psi^N({\mathcal I}_{w,\eta})*_{\bar s}{\mathrm R}\Psi^M({\mathcal A}_{\lambda,\eta})\ident
{\mathrm R}\Psi^P({\mathcal I}_{w,\eta}*_\eta{\mathcal A}_{\lambda,\eta})
}$$
\end{lemma}

\PROOF
According to a theorem of Beilinson-Bernstein (see the theorem 4.7 
in \cite{Illusie}) we have an
isomorphism of perverse sheaves
$${\mathrm R}\Psi^{M\times N}({\mathcal A}_{\lambda,\eta}\boxtimes_\eta{\mathcal I}_{w,\eta})
\ident{\mathrm R}\Psi^M({\mathcal A}_{\lambda,\eta})\boxtimes_{\bar s}{\mathrm R}\Psi^N({\mathcal I}_{w,\eta}).$$
This induces an isomorphism between the pull-backs
$$p_1^*{\mathrm R}\Psi^{M\times N}({\mathcal A}_{\lambda,\eta}\boxtimes_\eta{\mathcal I}_{w,\eta})
\ident p_1^*({\mathrm R}\Psi^M({\mathcal A}_{\lambda,\eta})\boxtimes_{\bar s}
{\mathrm R}\Psi^N({\mathcal I}_{w,\eta}))$$
which are up to the shift by the relative dimension 
$p_1$, perverse too. 
By definition, we have
$$p_1^*({\mathrm R}\Psi^{M}({\mathcal A}_{\lambda,\eta})
\boxtimes_{\bar s}{\mathrm R}\Psi^N({\mathcal I}_{w,\eta}))
\ident p_2^*({\mathrm R}\Psi^{M}({\mathcal A}_{\lambda,\eta})
\,{\tilde\boxtimes}_{\bar s}\,{\mathrm R}\Psi^N({\mathcal I}_{w,\eta})).$$
As $p_1, p_2$ are smooth, $p_1^*$ and $p_2^*$ commute with 
nearby cycle, so applying ${\mathrm R}\Psi^{{\tilde M} \times {\tilde N}}$ to
$$
p^*_1({\mathcal A}_{\lambda,\eta} \boxtimes_{\eta} {\mathcal I}_{w,\eta}) \ident p^*_2({\mathcal A}_{\lambda,\eta} {\tilde \boxtimes}_{\eta} {\mathcal I}_{w,\eta})
$$
gives an isomorphism 
$$p_1^*{\mathrm R}\Psi^{M\times N}({\mathcal A}_{\lambda,\eta}\boxtimes_\eta{\mathcal I}_{w,\eta})
\ident p_2^*{\mathrm R}\Psi^{M\,{\tilde\times}\,N}({\mathcal A}_{\lambda,\eta}\,
{\tilde\boxtimes}_\eta\,{\mathcal I}_{w,\eta}).$$
Since $p_2$ is smooth with connected geometric fibres,
the uniqueness part of theorem 4.2.5 of Beilinson-Bernstein-Deligne \cite{BBD} 
implies that we have an isomorphism
$${\mathrm R}\Psi^{M\,{\tilde\times}\,N}({\mathcal A}_{\lambda,\eta}\,{\tilde\boxtimes}_\eta
\,{\mathcal I}_{w,\eta})\ident {\mathrm R}\Psi^M({\mathcal A}_{\lambda,\eta})\,
{\tilde\boxtimes}_{\bar s}\,{\mathrm R}\Psi^N({\mathcal I}_{w,\eta}).$$

By applying now the functor ${\mathrm R} m_{*}$, we have an isomorphism
$${\mathrm R} m_{*}{\mathrm R}\Psi^{M\,{\tilde\times}\,N}({\mathcal A}_{\lambda,\eta}\,
{\tilde\boxtimes}_\eta\,{\mathcal I}_{w,\eta})
\ident{\mathrm R}\Psi^M({\mathcal A}_{\lambda,\eta})*_{\bar s}{\mathrm R}\Psi^N({\mathcal I}_{w,\eta}).$$
Since the functor ${\mathrm R}\Psi$ commutes with the direct image of a proper
morphism, we have
$${\mathrm R}\Psi^P({\mathcal A}_{\lambda,\eta}*_\eta{\mathcal I}_{w,\eta})\ident
{\mathrm R} m_{*}{\mathrm R}\Psi^{M\,{\tilde\times}\,N}({\mathcal A}_{\lambda,\eta}\,
{\tilde\boxtimes}_\eta\,{\mathcal I}_{w,\eta}).$$
By composing the above isomorphisms, we get
$${\mathrm R}\Psi^M({\mathcal A}_{\lambda,\eta})*_{\bar s} {\mathrm R}\Psi^N({\mathcal I}_{w,\eta})\ident
{\mathrm R}\Psi^P({\mathcal A}_{\lambda,\eta}*_\eta{\mathcal I}_{w,\eta}).$$
By the same argument, we prove
$${\mathrm R}\Psi^N({\mathcal I}_{w,\eta})*_{\bar s}{\mathrm R}\Psi^M({\mathcal A}_{\lambda,\eta})\ident
{\mathrm R}\Psi^P({\mathcal I}_{w,\eta}*_\eta{\mathcal A}_{\lambda,\eta}).$$
This finishes the proof of the lemma. $\square$

\bigskip
Now it clearly suffices to prove 
$$ {\mathcal A}_{\lambda,\eta}*_\eta{\mathcal I}_{w,\eta}\ident 
{\mathcal I}_{w,\eta}*_\eta{\mathcal A}_{\lambda,\eta}$$
which is an easy consequence of the following lemma.

\begin{lemma}
\begin{enumerate}
\item
Over the generic point $\eta$, we have two commutative 
triangles
$$\diagramme{
 & M_{\eta}\,{\tilde\times}\,N_{\eta} &  \cr
\hfill{}^{i}\!\swarrow &                  &\searrow\!^{m}\hfill \cr
M_{\eta}\,{\times}\,N_{\eta}  &
{\smash{\mathop{\hbox to 15mm{\rightarrowfill}}
     \limits^{\scriptstyle j}}}
&\,\ P_{\eta} \cr
\hfill_{i'}\!\nwarrow & &\nearrow\!_{ m'}\hfill \cr
 & N_{\eta}\,{\tilde\times}\,M_{\eta} &  \cr
}$$
where all arrows are isomorphisms.
\item
Morever, we have the following isomorphisms
$$\displaylines{
i^*({\mathcal A}_{\lambda,\eta}\boxtimes{\mathcal I}_{w,\eta})\ident
{\mathcal A}_{\lambda,\eta}\,{\tilde\boxtimes}\,{\mathcal I}_{w,\eta}\cr
i'{}^*({\mathcal A}_{\lambda,\eta}\boxtimes {\mathcal I}_{w,\eta}) \ident
{\mathcal I}_{w,\eta}\,{\tilde \boxtimes}\,{\mathcal A}_{\lambda,\eta}\cr}$$. 
\end{enumerate}
\end{lemma}

\subsection{Proof of Lemma 23}

Let us prove the above lemma in the linear case.

Over the generic point $\eta$, we have the canonical decomposition of
$${\bar{\mathcal V}}_{F}=t^{n_-}(t+\varpi)^{n'_--1}F[t]^d
/t^{n_+}(t+\varpi)^{n'_+}F[t]^d$$
into the direct sum
${\bar{\mathcal V}}_{F}={\bar{\mathcal V}}^{\,(t)}_{F}
\oplus {\bar{\mathcal V}}^{\,(t+\varpi)}_{F}$
where 
$$\displaylines{
{\bar {\mathcal V}}^{\,(t)}_{F}=t^{n_-}F[t]^d/t^{n_+}F[t]^d\cr
{\bar {\mathcal V}}^{\,(t+\varpi)}_{F}=
(t+\varpi)^{n'_- -1}F[t]^d/(t+\varpi)^{n'_+}F[t]^d.}$$
With respect to this decomposition, all the terms of the filtration
$${\bar {\mathcal V}}_0\subset{\bar {\mathcal V}}_1\subset\cdots\subset
{\bar{\mathcal V}}_{d-1}$$
decompose to
${\bar{\mathcal V}}_i=
{\bar{\mathcal V}}^{\,(t)}_i\oplus{\bar{\mathcal V}}^{\,(t+\varpi)}_i$
for all $i=0,\ldots,d-1$. Here, we have
$${\bar{\mathcal V}}^{\,(t)}_0=\cdots={\bar{\mathcal V}}^{\,(t)}_{d-1}
=F[t]^d/t^{n_+}F[t]^d.$$
Let $R$ be an $F$-algebra and let $({\mathcal L}_\bullet,{\mathcal L}'_\bullet)$
be an element of $(M\,{\tilde\times}\,N)(R)$.
These chains of $R[t]$-modules verify
$$\displaylines{
t^{n_+}{\mathcal V}_{i,R}\subset{\mathcal L}_i\subset t^{n_-}{\mathcal V}_{i,R}\cr
(t+\varpi)^{n'_+}{\mathcal L}_{i}\subset{\mathcal L}'_i\subset (t+\varpi)^{n'_-}{\mathcal L}_{i}}$$
As usual, let ${\bar{\mathcal L}}_i,{\bar{\mathcal L}}'_i$ denote the image of ${\mathcal L}_i,{\mathcal L}'_i$ 
in ${\bar {\mathcal V}}_{R}$. As $R[t]$-modules, they decompose to
${\bar{\mathcal L}}_i={\bar{\mathcal L}}^{\,(t)}_i\oplus{\bar{\mathcal L}}^{\,(t+\varpi)}_i$  and 
${\bar{\mathcal L}}'_i={\bar{\mathcal L}}'{}^{\,(t)}_i\oplus{\bar{\mathcal L}}'{}^{\,(t+\varpi)}_i.$
The above inclusion conditions imply indeed 
$$
{\bar{\mathcal L}}^{\,(t)}_i={\bar{\mathcal L}}'{}^{\,(t)}_i\ ;\ 
{\bar{\mathcal L}}^{\,(t+\varpi)}_i={\bar{\mathcal V}}^{\,(t+\varpi)}_{i,R}.$$

Consequently, ${\mathcal L}_\bullet$ is completely determined by ${\mathcal L}'_\bullet$.
In other terms, the map 
$m({\bar{\mathcal L}}_\bullet,{\bar{\mathcal L}}'_\bullet)={\bar{\mathcal L}}'_\bullet$
is an isomorphism of functors over $\eta$. 
In the same way, the map
$$i({\bar{\mathcal L}}_\bullet^{\,(t)}\oplus{\bar{\mathcal V}}^{\,(t+\varpi)}_{\bullet,R},
{\bar{\mathcal L}}_\bullet^{\,(t)}\oplus{\bar{\mathcal L}}'{}^{\,(t+\varpi)}_\bullet)=
({\bar{\mathcal L}}_\bullet^{\,(t)}\oplus{\bar{\mathcal V}}^{\,(t+\varpi)}_{\bullet,R},
{\bar{\mathcal V}}_{\bullet,R}^{\,(t)}\oplus{\bar{\mathcal L}}'{}^{\,(t+\varpi)}_\bullet)$$
yields an isomorphism 
$i:M_{\eta}\,{\tilde\times}\, N_{\eta}
\ident M_{\eta}{\times} N_{\eta}$.
The composed isomorphism $j=m\circ i^{-1}$ is given by
$$j({\bar{\mathcal L}}_\bullet^{\,(t)}\oplus{\bar{\mathcal V}}^{\,(t+\varpi)}_{\bullet,R},
{\bar{\mathcal V}}_{\bullet,R}^{\,(t)}\oplus{\bar{\mathcal L}}'{}^{\,(t+\varpi)}_\bullet)
={\bar{\mathcal L}}{}^{\,(t)}_\bullet\oplus
{\bar{\mathcal L}}'{}^{\,(t+\varpi)}_\bullet.$$
The analogous statement for
the lower triangle in the diagram can be proved in the same way
and the first part of the lemma is proved.

\bigskip

By the very definition of 
${\mathcal A}_{\lambda,\eta}\,{\tilde\boxtimes}\,{\mathcal I}_{w,\eta}$,
in order to prove the second part of the lemma,
it suffices to construct an isomorphism
$$p_1^*({\mathcal A}_{\lambda,\eta}\boxtimes{\mathcal I}_{w,\eta})\ident
p_2^*\,i^*({\mathcal A}_{\lambda,\eta}\boxtimes{\mathcal I}_{w,\eta}).$$
In fact, the triangle
$$\diagramme{
 & {\tilde M}_{\eta}\times {\tilde N}_{\eta} &  \cr
\hfill{}^{p_1}\!\swarrow &  &\searrow\!^{p_2}\hfill \cr
M_{\eta}\,{\times}\,N_{\eta}     &
{\smash{\mathop{\hbox to 15mm{\leftarrowfill}}
                               \limits^{\scriptstyle i}}} &
 M_{\eta}\,{\tilde\times}\,N_{\eta}}$$
does not commute. Nevertheless this lack of commutativity
can be corrected by equivariant properties.
We consider the diagram
$$\diagramme{
&{\tilde M}_{\eta}\times{\tilde N}_{\eta}&&\cr
\hfill{}^{q_1}\!\swarrow && \searrow\!^{q_2}\hfill     &\cr
\hfill {\tilde J}_{\eta}\times M_{\eta}\,{\times}\,N_{\eta}\!\!\!\!\!\! &
{\smash{\mathop{\hbox to 10mm{\leftarrowfill}}
\limits^{\scriptstyle{\mathrm Id}\times i}}} &
\!\!\!\!\!\!{\tilde J}_{\eta}\times M_{\eta}\,{\tilde\times}\,N_{\eta}&\cr
{}^{{\mathrm pr}_1}\!\swarrow & &\!\!\!\!\!\!\swarrow\!{}_{{\mathrm pr}_2}\hfill
\searrow \!{}^{\alpha}&\cr
M_{\eta}\,{\times}\,N_{\eta}\ \ \
{\smash{\mathop{\hbox to 15mm{\leftarrowfill}}
\limits^{\scriptstyle i}}}\!\!\!\! &
\ M_{\eta}\,{\tilde\times}\,N_{\eta}& &
\!\!\!\!
\!\!\!\!\!\!M_{\eta}\,{\tilde\times}\,N_{\eta}\ \ \cr
}$$
defined as follows.

For any $F$-algebra $R$, an element $g\in {\tilde M}(R)$
is an $R[t]$-endomorphism of ${\bar {\mathcal V}}_{R}$ such that 
${\bar{\mathcal L}}_\bullet=g(t^{n_-}{\bar{\mathcal V}}_{\bullet,R})\in M(R)$.
As ${\bar {\mathcal V}}_{R}$ decomposes to
${\bar {\mathcal V}}_{R}={\bar {\mathcal V}}_{R}^{\,(t)}\oplus
{\bar {\mathcal V}}_{R}^{\,(t+\varpi)}$,
its $R[t]$-endomorphism $g$ can be identified to a pair 
$g=(g^{\,(t)},g^{\,(t+\varpi)})$ 
where $g^{\,(t)}$, respectively 
$g^{\,(t+\varpi)}$, is an endomorphism of ${\bar {\mathcal V}}_{R}^{\,(t)}$,
respectively of ${\bar {\mathcal V}}_{R}^{\,(t+\varpi)}$.

As we have seen above, for $\bar{\mathcal L}_\bullet\in M(R)$, we have 
${\bar{\mathcal L}}_i= {\bar{\mathcal L}}_i^{\,(t)}\oplus {\bar{\mathcal L}}_i^{\,(t+\varpi)}$
with ${\bar{\mathcal L}}_i^{\,(t+\varpi)}={\bar {\mathcal V}}_{i,R}^{\,(t+\varpi)}$. 
Consequently, $g^{\,(t+\varpi)}$ is an automorphism of
${\bar {\mathcal V}}_{R}^{\,(t+\varpi)}$ fixing the filtration 
${\bar{\mathcal V}}_{\bullet,R}^{\,(t+\varpi)}$.
In a similar way, an element $g'\in {\tilde N}(R)$ can be identified
with a pair $(g'{}^{\,(t)},g'{}^{\,(t+\varpi)})$ 
where $g'{}^{\,(t)}$ is an automorphism of ${\bar {\mathcal V}}_{R}^{\,(t)}$ 
fixing the filtration ${\bar {\mathcal V}}^{(t)}_{\bullet,R}$.

\begin{itemize}
\item The morphism $q_1$ is defined by
$$q_1(g,g')=((g'{}^{\,(t)},g^{\,(t+\varpi)}),
g^{\,(t)}t^{n_-}{\bar{\mathcal V}}_{\bullet,R}^{\,(t)}\oplus{\bar{\mathcal V}}_{\bullet,R}^{\,(t+\varpi)},
{\bar{\mathcal V}}_{\bullet,R}^{\,(t)}\oplus 
g'{}^{\,(t+\varpi)}(t+\varpi)^{n'_-}{\bar{\mathcal V}}_{\bullet,R}^{\,(t+\varpi)}).$$

\item The morphism $q_2$ is defined by
$$q_2(g,g')=((g'{}^{\,(t)},g^{\,(t+\varpi)}),
g^{\,(t)}t^{n_-}{\bar{\mathcal V}}_{\bullet,R}^{\,(t)}\oplus{\bar{\mathcal V}}_{\bullet,R}^{\,(t+\varpi)},
g^{\,(t)}t^{n_-}{\bar{\mathcal V}}_{\bullet,R}^{\,(t)}\oplus 
g'{}^{\,(t+\varpi)}(t+\varpi)^{n'_-}{\bar{\mathcal V}}_{\bullet,R}^{\,(t+\varpi)}).$$

\item The morphism $\alpha$ is defined by
$$\displaylines{
\alpha((g'{}^{\,(t)},g^{\,(t+\varpi)}),
{\bar{\mathcal L}}_\bullet^{\,(t)}\oplus{\bar{\mathcal V}}_{\bullet,R}^{\,(t+\varpi)},
{\bar{\mathcal L}}_\bullet^{\,(t)}\oplus 
{\bar{\mathcal L}'}_\bullet{}^{\,(t+\varpi)})\cr
=({\bar{\mathcal L}}_\bullet^{\,(t)}\oplus{\bar{\mathcal V}}_{\bullet,R}^{\,(t+\varpi)},
{\bar{\mathcal L}}_\bullet^{\,(t)}\oplus 
g^{\,(t+\varpi)}{\bar{\mathcal L}'}_\bullet{}^{\,(t+\varpi)}).}$$

\item ${\mathrm pr}_1$ and ${\mathrm pr}_2$ are the obvious projections 
\end{itemize}

We can can easily check that this diagram commutes
and that
$${\mathrm pr}_1\circ q_1= p_1\ ;\ \alpha\circ q_2=p_2.$$
Now it is clear that
$$p_1^*({\mathcal A}_{\lambda,\eta}\boxtimes{\mathcal I}_{w,\eta})
\ident q_2^*\,{\mathrm pr}_2^*\,i^* 
({\mathcal A}_{\lambda,\eta}\boxtimes{\mathcal I}_{w,\eta}).$$
Moreover, by equivariant properties of ${\mathcal A}_\lambda$ and ${\mathcal I}_w$,
we have
$${\mathrm pr}_2^*\,i^* 
({\mathcal A}_{\lambda,\eta}\boxtimes{\mathcal I}_{w,\eta})
\ident \alpha^*i^* 
({\mathcal A}_{\lambda,\eta}\boxtimes{\mathcal I}_{w,\eta}).$$
(Note that the group $I_\eta$ acts on $M_\eta {\tilde \times} N_\eta$ by acting on the second factor of $M_\eta \times N_\eta \cong M_\eta {\tilde \times} N_\eta$ and $\alpha$ gives the corresponding action of ${\tilde J}_\eta$ via the projection ${\tilde J}_\eta \rightarrow I_\eta$.)
In putting these things together, we get the required isomorphism
$$p_1^*({\mathcal A}_{\lambda,\eta}\boxtimes{\mathcal I}_{w,\eta})
\ident p_2^*\,i^* 
({\mathcal A}_{\lambda,\eta}\boxtimes{\mathcal I}_{w,\eta}).$$
This finishes the proof of the lemma in the linear case.

In the symplectic case, let us mention that the $F$-vector  space
$$t^{n_-}(t+\varpi)^{n'_- -1}F[t]^{2d}/t^{n_+}(t+\varpi)^{n'_+}F[t]^{2d}$$
equipped with the symplectic form 
$t^{-n_- -n_+}(t+\varpi)^{-n'_- -n'_+ +1}\langle\,,\,\rangle $
splits into the direct sum of two vector spaces
$$t^{n_-}F[t]^{2d}/t^{n_+}F[t]^{2d}\oplus
(t+\varpi)^{n'_- -1}F[t]^{2d}/(t+\varpi)^{n'_+}F[t]^{2d}$$
equipped with symplectic forms
$t^{-n_--n_+}\langle\,,\,\rangle $
and $(t+\varpi)^{-n'_- -n'_+ +1}\langle\,,\,\rangle $ respectively.
Further, note that $g \in {\tilde M}(R)$ decomposes as $g = (g^{(t)},g^{(t+\varpi)})$ where
$g^{(t)} \in {\mathrm Aut}_{R[t]}(t^{n_-}R[t]^{2d} / t^{n_+}R[t]^{2d})$ is such that $\langle g^{(t)}x,g^{(t)}y \rangle = c_{g^{(t)}}t^{-n_- + n_+}\langle x,y \rangle$ (for some $c_{g^{(t)}} \in R^{\times}$), 
and $g^{(t+\varpi)} \in {\mathrm Aut}_{R[t]}((t+\varpi)^{n'_- -1}R[t]^{2d} / (t+\varpi)^{n'_+}R[t]^{2d})$ is such that $\langle g^{(t+\varpi)}x, g^{(t+\varpi)}y \rangle = c_{g^{(t+\varpi)}}\langle x,y \rangle$ (for some $c_{g^{(t+\varpi)}} \in R^{\times}$).  A similar decomposition 
$g' = (g'{}^{\,(t)},g'{}^{\,(t+\varpi)})$ holds, and thus ones sees $(g'{}^{\,(t)},g^{(t+\varpi)}) \in {\tilde J}(R)$.  Thus the maps $q_1$ and $q_2$ as defined above make sense in the symplectic case as well.  The rest of the argument goes through without change as in the linear case. 

This finishes the proof of Lemma 23.  We have therefore finished the proof of Proposition 21, and thus Proposition 13 and Theorem 11 as well.  $\square$

\bibliographystyle{plain}

\begin{thebibliography}{1}

\bibitem{BBD}
A.A.~Beilinson ; J.~Bernstein ; P.~Deligne. \newblock
{Faisceaux pervers} {\em in} 
{Analyse et topologie sur les espaces singuliers}, I  
\newblock{\em Ast\'risque} 100 (1982).

\bibitem{Dat}
J.-F. Dat.\newblock 
{ Caract\`eres \`a valeurs dans le centre de Bernstein}. 
\newblock{\em J. Reine Angew. Math.} 508 (1999), 61--83.

\bibitem{Deligne}
P.~Deligne.\newblock
{ Le formalisme des cycles \'evanescents} {\it in}
SGA 7 II, LNM 340, Springer 1973. 

\bibitem{Gaitsgory}
D.~Gaitsgory. \newblock {Construction of central elements in the affine Hecke algebra via nearby cycles}. {\em Preprint}. math.AG/9912074 9 Dec 1999.

\bibitem{Ginzburg}
V.~Ginzburg.\newblock { Perverse sheaves on a loop group and 
Langlands duality}.
\newblock {\em Preprint} (1996)

\bibitem{Grothendieck}
A.~Grothendieck.
{Formule de Lefschetz et rationalit\'e des fonctions $L$}
S\'em. Bourbaki no 279.



\bibitem{Haines}
T.~Haines. \newblock {The combinatorics of Bernstein functions.}
{\em  Preprint} (1998).

\bibitem{Haines2}
T.~Haines. \newblock {Test Functions for Shimura Varieties: the Drinfeld Case.} {\em Preprint} (1998).

\bibitem{Illusie} L.~Illusie.\newblock 
{ Autour du th\'eor\`eme de monodromie locale}
{\em in} { P\'eriodes $p$-adiques}. 
\newblock{\em Ast\'erisque} No. 223 (1994), 9--57. 


\bibitem{Kato}
S.-I.~Kato.  \newblock{ Spherical functions and a $q$-analogue of 
Kostant's weight multiplicity formula}.
\newblock{\em Invent. Math.} 66 (1982), no. 3, 461--468. 

\bibitem{Kora}
R.~Kottwitz; M.~Rapoport. \newblock{Minuscule Alcoves for $Gl_n$ and $GSp_{2n}$}.
\newblock{\em Preprint} (1998).


\bibitem{Lusztig}
G.~Lusztig.\newblock { Singularities, characters formulas and a $q$-analogue of
weight multiplicities} {\em in} {Analyse et topologie sur les espaces singuliers.}
\newblock {\em Ast\'erisque} 101-102 (1983), 200-229.

\bibitem{Lusztig1}
G.~Lusztig.\newblock 
{Cells in affine Weyl groups and tensor categories.}
\newblock{\em Advances in Math.} 129(1997),  85--98.


\bibitem{Mathieu}
O.~Mathieu.\newblock  
{ Formules de caract\`eres pour les alg\`ebres de Kac-Moody g\'en\'erales}. 
\newblock {\em Ast\'erisque } 159-160 (1988), 267 pp. 

%\bibitem{Matsumura}
%H.~Matsumura. \newblock{ Commutative Ring Theory}.
%\newblock{\em Cambridge Studies in Adv. Math.} 8, Cambridge Univ. Press (1986), 320 pp.


\bibitem{Mirkovic-Vilonen}
I.~Mirkovic, K.~Vilonen.\newblock { Perverse sheaves on loop grassmannians
and Langlands duality}.\newblock{\it Preprint} (1997).


\bibitem{Rapoport} 
M.~Rapoport.\newblock{ On the bad reduction of the Shimura varieties}
{\it in} { Automorphic forms, Shimura varieties and $L$-functions}
edited by L.Clozel and J. Milne.
\newblock{\em Perspec. Math.} 11, p. 253-321, Acad. Press 1990.

\bibitem{Rapoport-Zink}
M.~Rapoport,Th.~ Zink.\newblock { Period spaces for $p$-divisible groups}.
\newblock {\em Annals of Math. Studies} 144, Princeton Univ. Press 1996


\bibitem{Satake}
I.~Satake.
\newblock Theory of spherical functions on reductive algebraic 
groups over $p$-adic  fields.
\newblock {\em Publ. IHES}, 18: 1--69, 1963.




\end{thebibliography}

\end{document}